\documentclass[12pt]{amsart}
\usepackage{amssymb}
\usepackage{times}
\usepackage{graphicx}

\theoremstyle{plain}

\numberwithin{equation}{section}

\newcommand{\calB}{\mathcal{B}}

\newcommand{\calD}{\mathcal{D}}

\newcommand{\calH}{\mathcal{H}}

\newcommand{\calO}{\mathcal{O}}

\newcommand{\bbF}{\mathbb{F}}
\newcommand{\bbC}{\mathbb{C}}

\newcommand{\bbP}{\mathbb{P}}
\newcommand{\bbQ}{\mathbb{Q}}

\newcommand{\bbZ}{\mathbb{Z}}

\newcommand{\la}{\langle}
\newcommand{\ra}{\rangle}

\def\SL{{\text{SL}}}
\def\Aut{{\text{Aut}}}

\def\Ker{{\text{Ker}}}
\def\O{{\text{O}}}

\def\GL{{\text{GL}}}
\def\Sp{{\text{Sp}}}

\def\SL{{\text{SL}}}

\def\mod{{\text{mod}}}
\def\lim{{\text{lim}}}
\def\sin{{\text{sin}}}
\def\Hom{{\text{Hom}}}

\begin{document}
\title [Moduli of plane quartics] {Moduli of plane quartics, G\"opel invariants and Borcherds products}
\author{Shigeyuki Kond{$\bar{\rm o}$}}
\address{Graduate School of Mathematics, Nagoya University, Nagoya,
464-8602, Japan}
\email{kondo@math.nagoya-u.ac.jp}
\thanks{Research of the author is partially supported by
Grant-in-Aid for Scientific Research A:18204001, Houga:20654001, Japan}

\begin{abstract}
It is known that the moduli space of plane quartic curves is birational to an arithmetic
quotient of a 6-dimensional complex ball (\cite{Kon1}).  In this paper,
we shall show that there exists a 15-dimensional space of meromorphic automorphic forms on the complex ball which gives a birational embedding of the moduli space of plane quartics with level 2 structure into $\bbP^{14}$.  This map coincides with the one given by Coble \cite{C} by using G{\" o}pel invariants.

\end{abstract}
\maketitle

CONTENTS
\begin{itemize}
\item[1.] Introduction
\item[2.] Preliminaries
\item[3.] Plane quartics and $7$ points on $\bbP^2$
\item[4.] G{\" o}pel invariants
\item[5.] $K3$ surfaces associated to plane quartic curves
\item[6.] Hermitian form and reflections
\item[7.] Heegner divisors
\item[8.] Weil representation
\item[9.] Automorphic forms: Additive liftings
\item[10.] Automorphic forms: Borcherds products
\item[11.] Automorphic forms and G{\" o}pel invariants
\end{itemize}

\medskip

\section{Introduction}
In \cite{Kon1}, the author proved that the moduli space of smooth plane quartics is isomorphic to an arithmetic quotient 
$(\calB \setminus \calH)/\Gamma$ where $\calB$ is a 6-dimensional complex ball, $\calH$ is the union of
hypersurfaces in $\calB$ and $\Gamma$ is an arithmetic subgroup of $\Aut (\calB )$.  Moreover $\calH$ decomposes into two types: $\calH = \calH_n \cup \calH_h$ where 
a generic point in $\calH_n$ (resp. $\calH_h$) corresponds to a plane quartic with a node 
(resp. a hyperelliptic curve of genus 3).  The above isomorphism is defined as follows.  Let $C$ be a smooth plane quartic.  
Then taking the 4-cyclic cover of $\bbP^2$ branched
along $C$, we have a $K3$ surface $X$ and an automorphism $\sigma$ of $X$ of order 4.  
By the theory of periods of $K3$ surfaces, we can see that the period domain
of the pairs $(X, \sigma)$ is $\calB \setminus \calH$.

There exists a subgroup $\tilde{\Gamma}$ 
of $\Gamma$ such that $\Gamma/\tilde{\Gamma} \cong \bbZ/2\bbZ \cdot \Sp(6,\bbF_2) $ and
the quotient $(\calB \setminus \calH)/(\tilde{\Gamma} \cdot \bbZ/2\bbZ)$ is $\Sp(6,\bbF_2)$-equivariantly 
isomorphic to the moduli space of smooth plane
quartics with level 2-structure (Proposition \ref{kondo}).  On the other hand, a smooth plane quartic $C$ is
naturally corresponding to a del Pezzo surface $S$ of degree two by taking the double cover
of $\bbP^2$ branched along $C$.  Note that $S$ is nothing but the quotient of the above $X$ by the involution $\sigma^2$.
Conversely the anti canonical model of $S$ is the double cover of $\bbP^2$ branched along a smooth quartic.  Thus the moduli space of smooth
plane quartics is isomorphic to the moduli of del Pezzo surface of degree two.  
A del Pezzo surface $S$ of degree two is obtained by blowing ups at seven points on $\bbP^2$
in general position.  We call the contraction of $S\to \bbP^2$ a marking.  Here we consider an order of 7 points and
7 exceptional curves.  It is known that  there are $2^{10}\cdot 3^4\cdot 5\cdot 7$ markings.  
The group of changing markings is isomorphic to $\Sp(6,\bbF_2) \cong W(E_7)/\{\pm 1\}$ where $W(E_7)$ is the Weyl group of type $E_7$.
This implies that the moduli space of smooth plane quartics with level 2-structure 
is $\Sp(6,\bbF_2)$-equivariantly isomorphic to the moduli of ordered 7 points on $\bbP^2$ 
in general position.  Let $P^7_2$ be the moduli space of semi-stable ordered 7 points on 
$\bbP^2$.
Thus we have a birational, $\Sp(6,\bbF_2)$-equivariant isomorphism between $P^7_2$ and 
$(\calB \setminus \calH)/(\tilde{\Gamma} \cdot \bbZ/2\bbZ)$ (Proposition \ref{kondo1}).  

In \cite{C}, by using G\"opel invariants, Coble showed that there exists a $\Sp(6,\bbF_2)$-equivariant birational embedding of $P^7_2$ into $\bbP^{14}$ whose image satisfies 63
cubic relations.

In this paper, by using the theory of automorphic forms due to Borcherds \cite{B1}, \cite{B2}, Freitag \cite{F}, 
we shall show that there exists a 15-dimensional space of {\it meromorphic} automorphic forms on $\calB$ which gives a $\Sp(6,\bbF_2)$-equivariant birational embedding of $\calB/(\tilde{\Gamma} \cdot \bbZ/2\bbZ)$
into $\bbP^{14}$ (Theorem \ref{main4}).  This map coincides with the one given by Coble \cite{C}.

The plan of this paper is as follows.  In \S 2 we fix the notation of lattices.  In \S 3,4, we recall classical results on the moduli of
7 points on $\bbP^1$ and Coble's work \cite{C}, \cite{DO}.  
In \S 5,6,7 we study periods of the pairs $(X, \sigma)$.  The sections 8 and 9 devote to construct a linear system of automorphic forms of dimension 15 (Theorem \ref{main1}).  
To determine the divisor
of a member of this linear system, we need a very special automophic form with known zeros and poles.  In \S 10 we show the existence
of such an automorphic form (Theorem \ref{multiplicative}, Corollary \ref{multi}, Theorem \ref{main2}, Corollary \ref{main3}).  
Finally in \S 11, we shall give the main theorem (Theorem \ref{main4}) and its proof.

We mention the related works.  Main idea in this paper follows from the paper Allcock-Freitag \cite{AF} in which they considered the case
of del Pezzo surfaces of degree 3, that is, cubic surfaces.   They gave a $W(E_6)$-equivariant embedding of the moduli of marked cubic surfaces into $\bbP^9$ by using automorphic forms.
In this case, roughly speaking, automorphic forms correspond to Cayley's cross ratios of cubic surfaces.   
In case of hyperelliptic curves of genus 3 with level 2-structure 
or equivalently ordered 8 points on $\bbP^1$, the author \cite{Kon3} gave a $S_8$-equivariant  
embedding of the moduli space into $\bbP^{13}$ by using Borcherds theory.  
In this case, automorphic forms correspond to cross ratios of ordered 8 points on $\bbP^1$.
Also the author \cite{Kon2} gave an $\O(10,\bbF_2)^+$-equivariant 
embedding of the moduli space of Enriques surfaces with "level 2-structure" into 
$\bbP^{185}$.  A mistake in \cite{Kon2} was corrected in \cite{FS}.  
The author does not know the geometric meaning of automorphic forms in case of
Enriques surfaces.  All of these cases used holomorphic automorphic forms.  On the other hand,
in this paper, the author uses meromorphic automorphic forms.

\medskip
\emph{Acknowledgments}:
The author thanks Igor Dolgachev, Eberhard Freitag and Riccardo Salvati-Manni for valuable conversations.

\medskip

\section{Preliminaries}\label{}

A {\it lattice} $(L, \la, \ra)$ is a pair of a free $\bbZ$-module $L$ of rank $r$ and
a non-degenerate symmetric bilinear form $\la, \ra $ on $L$. 
For simplicity we omit $\la, \ra$ if there are no confusion.  For $x \in L$, we call $x^2 =\la x,x\ra$ the {\it norm} of $x$.
For a lattice $(L,\la,\ra)$ and an integer
$m$, we denote by $L(m)$ the lattice $(L, m\la, \ra)$.  We denote by $U$ the lattice
$(\bbZ^{\oplus 2}, 
\begin{pmatrix}0&1
\\1&0
\end{pmatrix})$ and by $A_m, \ D_n$ or $\ E_k$ the even {\it negative} definite lattice defined by
the Cartan matrix of type $A_m, \ D_n$ or $\ E_k$ respectively.  For an integer $m$, we denote by $\la m\ra$ the lattice of rank 1
generated by a vector with norm $m$.  We denote by $L\oplus M$ the orthogonal direct sum of lattices $L$ and $M$.
We also denote by $L^{\oplus k}$ the orthogonal direct sum of $k$ copies of $L$.

Let $L$ be an even lattice and let $L^* =\Hom(L,\bbZ)$.  We denote by $A_L$ the quotient
$L^*/L$ and define a map
$$q_L : A_L \to \bbQ/2\bbZ, \ b_L : A_L\times A_L \to \bbQ/\bbZ$$
by $q_L(x+L) = \la x, x\ra\ \mod\ 2\bbZ$ and $b_L(x+L,y+L) =
\la x, y\ra\ \mod \ \bbZ$.  We call $q_L$, $b_L$ the {\it discriminant quadratic form, discriminant
bilinear form}, respectively.

Let $\O(L)$ be the orthogonal group of $L$, that is, the group of isomorphisms of $L$ preserving the bilinear form.
Similarly $\O(q_L)$ denotes the group of isomorphisms of $A_L$ preserving $q_L$.
There is a natural map
$$\O(L) \to \O(q_L)$$
whose kernel is denoted by $\tilde{\O}(L)$.
For more details we refer the reader to \cite{N1}.

\medskip

\section{Plane quartics and 7 points on $\bbP^2$}

Let $C$ be a smooth plane quartic curve in $\bbP^2$.  Then $C$ is a smooth curve of genus three.  
Conversely the canonical model of a general smooth curve of genus three is a plane quartic.  
An algebraic surface $S$ is called a {\it del Pezzo surface of degree} $d$ 
if the anti canonical class $-K_S$ is ample and $(-K_S)^2=d$.
A smooth plane quartic $C$ is
naturally corresponding to a del Pezzo surface $S$ of degree two by taking the double cover
of $\bbP^2$ branched along $C$.  Conversely the anti canonical model of $S$ is the double cover of $\bbP^2$ branched along a smooth quartic.  Thus the moduli space of smooth
plane quartics is isomorphic to the moduli of del Pezzo surface of degree two.
It is well known that there are 28 bitangent lines to $C$.  On $S$, bitangent lines
split into 56 $(-1)$-curves.

A del Pezzo surface $S$ of degree two is obtained by blowing ups at seven points on $\bbP^2$
in general position.  Here seven points $p_1,..., p_7$ in $\bbP^2$ are {\it in general position} if
no two points coincide, no three points lie on a line and no six points lie on a conic.
Let $e_1,..., e_7$ be exceptional curves and $e_0$ the total transform of the line.
The Picard lattice of $S$ is isomorphic to $\la 1\ra \oplus \la -1\ra^{\oplus 7}$ where $\la 1\ra$ is generated by $e_0$
and $\la -1\ra^{\oplus 7}$ is generated by $e_1,..., e_7$.  The orthogonal complement 
of the anti-canonical class $3e_0 - e_1 - \cdot \cdot \cdot - e_7$
is isomorphic to the root lattice $E_7$. 
There are 56 $(-1)$-curves, that is, 7 exceptional curves, proper transforms of 21 lines through
two points from $p_1,...,p_7$, proper transforms of 21 conics throuh 5 points from
$p_1,...,p_7$, proper transforms of 7 cubics through 7 points and having a node at one of
7 points:
\begin{equation}\label{discriminant condition}
\begin{cases}
e_i  \ (1\leq i \leq 7),\cr
e_0-e_i-e_j \ (1\leq i<j\leq 7),\cr
2e_0-e_1-e_2-e_3-e_4-e_5-e_6-e_7 + e_i+e_j \ (1\leq i<j \leq 7),\cr
3e_0-e_1-e_2-e_3-e_4-e_5-e_6-e_7 - e_i \ (1\leq i \leq 7).\cr
\end{cases}
\end{equation}
On the other hand, $p_1,..., p_7$ are in general position iff all of the 63 divisors
\begin{equation}\label{discriminant condition}
\begin{cases}
e_i - e_j \ (1\leq i < j \leq 7),\cr
e_0-e_i-e_j-e_k \ (1\leq i<j<k \leq 7),\cr
2e_0-e_1-e_2-e_3-e_4-e_5-e_6-e_7 + e_i \ (1\leq i\leq 7)\cr
\end{cases}
\end{equation}
on $S$ are not effective.  We call (\ref{discriminant condition}) {\it the discriminant conditions} of 7 points in $\bbP^2$.
Note that 63 divisors in the discriminant condition (\ref{discriminant condition}) 
can be considered as a set of positive roots of $E_7$.
There are $2^{10}\cdot 3^4\cdot 5\cdot 7$ sets of disjoint seven $(-1)$-curves on $S$, 
in other words, there are 
$2^{10}\cdot 3^4\cdot 5\cdot 7$
contractions of $S$ to $\bbP^2$.  We call a contraction $S \to \bbP^2$ a {\it marking} of $S$.  The group of
changing markings is isomorphic to $\Sp(6,\bbF_2)$.  Let $W(E_7)$ be the Weyl group of $E_7$.  
Then there is an exact sequence
\begin{equation}\label{}
1 \to \{\pm 1\} \to W(E_7) \to \Sp(6,\bbF_2)\to 1
\end{equation}
where $\{\pm 1\}$ is generated by the covering transformation of $S\to \bbP^2$.
A marking of $S$ naturally corresponds 
to a level 2 structure of a smooth curve  $C$ of genus 3.  Thus the moduli space of
del Pezzo surfaces of degree 2 with a marking is isomorphic to the moduli of smooth plane quartics with a level 2 structure.  This isomorphism is $\Sp(6,\bbF_2)$-equivariant.
Let $P_2^7$ be the moduli space of semi-stable ordered 7 points on $\bbP^2$. 
Here a set of ordered 7 points $\{ p_1,..., p_7\}$ is {\it semi-stable} (={\it stable}) iff
at most two points coincide and at most 4 points are collinear (e.g. \cite{DO}, p.120).
The discriminat condition (\ref{discriminant condition}) gives 63 subvarieties in $P_2^7$
whose complement parametrizes smooth plane quartics.  

For more details, we refer the reader to \cite{DO}.

\medskip

\section{G\"opel invariants}

In the following we shall recall Coble's work \cite{C} (also see \cite{DO}) on G\"opel invariants.  
For simplicity we denote by $\alpha (i,j)$, $\alpha (i,j,k)$, $\alpha (i)$ the positive roots
$e_i - e_j$, $e_0-e_i-e_j-e_k$, $2e_0-e_1-e_2-e_3-e_4-e_5-e_6-e_7 + e_i$, respectively.

The structure of lattice on $E_7$ induces a structure of quadratic form on $E_7/2E_7$ over
$\bbF_2$ and of dimension 7.  It has a radical $\sum_{i=0}^7 e_i \ \mod\ 2E_7$ and modulo radical it induces a
 symplectic space $\bbF_2^6$ of dimension 6.

A maximal totally isotropic subspace of $\bbF_2^6$ is called {\it G\"opel subspace}.
There are 135 G\"opel subspaces.  Each non-zero element in $\bbF_2^6$
corresponds to a positive root in $E_7$ and 7 non-zero elements in a G\"opel subspace gives
a set of mutually orthogonal 7 positive roots (\cite{DO}, Chap. IX, Lemma 8).
We call a set of mutually orthogonal 7 positive roots a {\it G\"opel subset}.
There are 135 G\"opel subsets (\cite{C}, \S 28):
\begin{equation}
\begin{cases}
30\ {\rm of\ type}\ \{ \alpha (1,2,3), \alpha(1,4,5), \alpha(2,4,6), \alpha(3,5,6), \alpha(1,6,7), \alpha(2,5,7), \alpha(3,4,7) \},\cr
105\ {\rm of\ type}\ \{ \alpha(1), \alpha(1,2,3), \alpha(1,4,5), \alpha(1,6,7), \alpha(2,3), \alpha(4,5),
\alpha(6,7) \}.\cr
\end{cases}
\end{equation}
For each G\"opel subset (or corresponding G\"opel subspace) $M$, 
we define a function $G_M$ on $(\bbC^3)^7$ called {\it G\"opel invariant} as follows.
First consider a G\"opel subset $M$ given by
\begin{equation}\label{goepel-1}
\{ \alpha (1,2,3), \alpha(1,4,5), \alpha(2,4,6), \alpha(3,5,6), \alpha(1,6,7), \alpha(2,5,7), \alpha(3,4,7) \}.
\end{equation}
Let $p_i, p_j, p_k$ be three points in $\bbP^2$ and let 
$v_i, v_j, v_k$ the column vectors in $\bbC^3$ corresponding to $p_i, p_j, p_k$ respectively.
We denote by $(i  j  k)$ the discriminant of the $3 \times 3$ matrix $(v_i v_j v_k)$.
Then the G\"opel invariant $G_M$ is defined by
\begin{equation}\label{goepel-ex}
G_M = (1 2 3)(1 4 5)(2 4 6)(3 5 6)(1 6 7)(2 5 7)(3 4 7).
\end{equation}
Obviously $G_M$ vanishes along 7 divisors corresponding to 7 roots in $M$.
Next consider the case $M$ is given by
\begin{equation}\label{goepel-2}
\{ \alpha (1), \alpha(1,2,3), \alpha(1,4,5), \alpha(1,6,7), \alpha(2,3), \alpha(4,5), \alpha(6,7)  \}.
\end{equation}
Then there are exactly three G\"opel subsets $M, M_1, M_2$ containing $\alpha(1,2,3), \alpha(1,4,5), \alpha(1,6,7)$:
\begin{equation}\label{goepel-3}
M_1= \{ \alpha (1,2,3), \alpha(1,4,5), \alpha(2,4,6), \alpha(3,5,6), \alpha(1,6,7), \alpha(2,5,7), \alpha(3,4,7) \}.
\end{equation}
\begin{equation}\label{goepel-4}
M_2= \{ \alpha (1,2,3), \alpha(1,4,5), \alpha(2,5,6), \alpha(3,5,7), \alpha(1,6,7), \alpha(2,4,7), \alpha(3,4,6) \}.
\end{equation}
Then $G_M$ is defined by
$$G_M= \pm(G_{M_1} - G_{M_2}) = (1 2 3)(1 4 5)(1 6 7)\{ (2 4 6)(3 5 6)(2 5 7)(3 4 7) - (2 5 6)(3 5 7)(2 4 7)(346)\}$$
up to sign $\pm$.
Then we can easily see that the second factor of $G_M$ vanishes along 4 subvarieties corresponding to
$ \alpha (1), \alpha(2,3), \alpha(4,5), \alpha(6,7)$.  For example, if we fix 5 points $p_2, p_3, p_4, p_5, p_6$ and
consider $p_7$ as a parameter, then the second factor of $G_M$ is a conic through these 5 points.

\subsection{Proposition}\label{Goepel}(Coble \cite{C}, \cite{DO}, Chap. IX)
(1) {\it Assume that $M$ is one of $30$ types, for example, $M$ is given by $(\ref{goepel-1})$.  Then $G_M$ vanishes exactly 
along $7$ divisors corresponding to $7$ roots in $M$ with multiplicity one.  Moreover $G_M$ vanishes along $21$ subvarieties corresponding to $21$ roots $\alpha (i,j)$ with multiplicity one.}

(2) {\it Assume that $M$ is one of $105$ types, for example, $M$ is given by $(\ref{goepel-2})$.  Then $G_M$ vanishes  
along $7$ subvarieties corresponding to $7$ roots in $M$ plus $21$ subvarieties corresponding to
$21$ roots $\alpha (i,j)$.  }

(3) {\it Let $A$ be a $2$-dimensional totally isotropic subspace.  Then there are exactly three G\"opel subspaces 
$M_1, M_2, M_3$ containing $A$.  Moreover $G_{M_i}\ (i=1,2,3)$ satisfy a linear relation }
$$G_{M_1}\pm G_{M_2} \pm G_{M_3} =0$$
where the sign $\pm$ is taken for a suitable one.
The {\it $135$ G\"opel invariants generate a $15$-dimensional space on which $\Sp(6,\bbF_2)$ acts linearly via its action on 
G\"opel subspaces.   This $15$-dimensional representation of $\Sp(6,\bbF_2)$ is irreducible.   }

\subsection{Theorem}\label{Coble}(Coble \cite{C}, Chap. IV, \cite{DO}, Chap. IX, Theorem 5)
{\it The $15$-dimensional linear system defines a $\Sp(6,\bbF_2)$-equivariant 
 rational map $\psi$ from $P_2^7$ to $\bbP^{14}$
which is birational onto its image.  The image satisfies $63$ cubic relations corresponding to
$63$ positive roots.}

\subsection{Lemma}\label{28/36}
{\it Each $G_M$ vanishes with multiplicity $1$ along exactly $28$ subvarieties among $36$ subvarieties defined by}:
\begin{equation}\label{4points}
\begin{cases}
e_0-e_i-e_j-e_k -e_l,  \ 1\leq i<j<k<l \leq 7,\cr
2e_0-e_1-e_2-e_3-e_4-e_5-e_6-e_7,\cr
\end{cases}
\end{equation}
{\it i.e. four points are collinear or seven points lie on a conic.}
\begin{proof}
Since $\Sp (6, \bbF_2)$ acts transitively on the set of G\"opel invariants, it is enough to see the assertion
for a G\"opel invariant.  Consider the case
$$G_M = (1 2 3)(1 4 5)(2 4 6)(3 5 6)(1 6 7)(2 5 7)(3 4 7).$$
Then $G_M$ does not vanish identically along 8 subvarieties defined by
$$2e_0-e_1-e_2-e_3-e_4-e_5-e_6-e_7, e_0-e_i-e_j-e_k -e_l$$ 
where
$$(i,j,k,l) = (1,2,4,7), (1,2,5,6), (1,3,4,6), (1,3,5,7), (2,3,4,5), (2,3,6,7), (4,5,6,7).$$
And $G_M$ vanishes along another 28 subvarieties with multiplicity 1.
\end{proof}
\medskip

Let $\pi : \hat{P}_2^7 \to P_2^7$ be the blow ups of $P_2^7$ along  
$21$ subvarieties corresponding to 21 roots $e_i - e_j$ ($1 \leq i < j \leq 7$) (\cite{DO}, Chap. IV). 
Thus $\hat{P}_2^7$ has 63 divisors corresponding to 63 discriminant conditions.  
For a positive root $\alpha$, we denote by $D_{\alpha}$ the corresponding divisor.
For G\"opel subsets $M_1, M_2$, we denote by $\hat{G}_{M_1}/\hat{G}_{M_2}$ the pull back of $G_{M_1}/G_{M_2}$ by $\pi$.  Then  

\subsection{Lemma}\label{blowup}
{\it The divisor of $\hat{G}_{M_1}/\hat{G}_{M_2}$ is given by}
$$\sum_{\alpha \in M_1} D_{\alpha} - \sum_{\beta \in M_2} D_{\beta}.$$

\begin{proof}
We remark that both $G_{M_1}$ and $G_{M_2}$ vanish along 21 subvarieties corresponding to $21$ roots $\alpha (i,j)$ 
with multiplicity one.  Hence the assertion follows from Proposition \ref{Goepel}.
\end{proof}

\medskip

Later we shall construct automorphic forms with the same property as G\"opel invariants.  See Corollaries
\ref{cubicrelation}, \ref{main3}

\medskip

\section{$K3$ surfaces associated to plane quartic curves}

\subsection{$K3$ surfaces and a bounded symmetric domain of type $IV$}

In \cite{Kon1}, the author proved that the moduli space of plane quartic curves is
birational to an arithmetic quotient of a 6-dimentional complex ball.  We recall this fact briefly.
Let $C$ be a smooth plane quartic curve and let $f_4(x,y,z) =0$ be the defining equation 
of $C$.  Consider the smooth quartic surface
$$X : f_4(x,y,z) + t^4 =0.$$
Obviously $X$ is a 4 cyclic cover of $\bbP^2$ branched along $C$ and is a $K3$ surface, that is, the canonical class
$K_X$ is trivial and $H^1(X, \calO_X) =0$.
If we denote by $S$ the double cover of $\bbP^2$ branched along $C$, then the 4 cyclic cover
factorizes and $X$ is the double cover of $S$ branched along $C$.
Let $\sigma$ be an automorphism of $X$ of order 4 which is the covering transformation
of $X\to \bbP^2$ and let $\iota = \sigma^2$.  The second cohomology group 
$H^2(X,\bbZ)$ together with the cup product adimits a structure of lattice. 
 As a lattice, $H^2(X,\bbZ)$ is an even unimodular lattice of signature $(3,19)$ which is uniquely determined by this property.
We denote by $L$ an abstract even unimodular lattice of signature $(3,19)$.  Then $H^2(X,\bbZ)\cong L$.
We define 
$$H^2(X,\bbZ)_+ = \{ x \in H^2(X,\bbZ) : \iota^*(x) =x\},\
H^2(X,\bbZ)_-=\{x \in H^2(X,\bbZ) : \iota^*(x) =-x\}.$$
Both $H^2(X,\bbZ)_{\pm}$ are 2-elementary lattices which are isomorphic to $L_{\pm}$
respectively where
\begin{equation}\label{basis}
L_+ = A_1(-1)\oplus A_1^{\oplus 7}, \ L_- = U\oplus U(2) \oplus D_4^{\oplus 2} \oplus A_1^{\oplus 2}.
\end{equation}
We denote by $\tilde{e}_i \ (0\leq i \leq 7)$ the pullback of $e_i$ by the map $X\to S$.
Then $\tilde{e}_0$ generates $A_1(-1)$ in (\ref{basis}) and $\tilde{e}_1,..., \tilde{e}_7$ generate $A_1^{\oplus 7}$.
Let
$$q_{L_+} : A_{L_+} = L_+^*/L_+ \to \bbQ /2\bbZ$$
be the discriminant quadratic form of $L_+$.  
Define $\tilde{\kappa} = 3\tilde{e}_0 - \tilde{e}_1 - \cdot \cdot \cdot - \tilde{e}_7$.
Then

\subsection{Lemma}\label{e7-1}(e.g. \cite{Kon1}, Lemma 2.1)  (i) {\it A vector $x$ in
$A_{L_+}$ with $q_{L_{+}}(x) = -1/2$ is represented by one of the following $56$ 
vectors}:
$$\tilde{e}_{i}/2, 
(\tilde{\kappa} - \tilde{e}_{i})/2, 1 \leq i \leq 7, 
(\tilde{e}_{0} - \tilde{e}_{i} - \tilde{e}_{j})/2, 
(\tilde{\kappa} - 
\tilde{e}_{0} + \tilde{e}_{i} + 
\tilde{e}_{j})/2, 1 \leq i < j \leq 7.$$
(ii)  {\it A vector $x$ in 
$A_{L_+}$ with $q_{L_{+}}(x) = 1/2$ is represented by one of the following $72$ 
vectors}:
$$\tilde{e}_{0}/2, (2\tilde{e}_0-\sum_{i} \tilde{e}_{i})/2, 
(\tilde{e}_{0}-\tilde{e}_i-\tilde{e}_j-\tilde{e}_k -\tilde{e}_l)/2,
(2\tilde{e}_{0}-\tilde{e}_{i} - \tilde{e}_{j} - \tilde{e}_{k})/2, 1 \leq i < j < k < l \leq 7.$$
(iii)  {\it A vector $x$ in $A_{L_+}$  with
$q_{L_{+}}(x) = 1$  is represented by one of the following
$64$  vectors}:
$${\tilde \kappa}/2, 
(\tilde{e}_{i} - \tilde{e}_{j})/2, 1 \leq i < j \leq 7, 
(\tilde{e}_{0} - \tilde{e}_{i} - \tilde{e}_{j} - \tilde{e}_{k})/2, 
1 \leq i < j < k \leq 7,$$
$$(2 \tilde{e}_{0} - \sum_{i} \tilde{e}_{i} + \tilde{e}_{j})/2,
1 \leq j \leq 7.$$
(iv)  {\it A vector $x$ in  $A_{L_+}$ with
$q_{L_{+}}(x) = 0$ is represented by one of the following
$64$  vectors}:
$$0, (\tilde{\kappa} + \tilde{e}_{i} - \tilde{e}_{j})/2, 1 \leq i < j
\leq 7,  
(\tilde{\kappa} + \tilde{e}_{0} - \tilde{e}_{i} - \tilde{e}_{j} - 
\tilde{e}_{k})/2,  1 \leq i < j < k \leq 7,$$
$$(\tilde{\kappa} + 2 \tilde{e}_{0} - \sum_{i} \tilde{e}_{i} +
\tilde{e}_{j})/2, 1 \leq j \leq 7.$$
(v)  {\it $\sigma^*$ acts trivially on vectors in {\rm (iii), (iv)} and acts on vectors $x$ in {\rm (i), (ii)} as}
$\sigma^*(x) = \tilde{\kappa}/2 - x$.

\subsection{Remark}\label{e7-3}
(1) Note that 63 vectors in Lemma \ref{e7-1}, (iii), except $\tilde{\kappa}/2$ correspond to
63 discriminant conditions (\ref{discriminant condition}).  Also 72 vectors in (ii) modulo the action of $\sigma^*$ 
correspond to 36 conditions (\ref{4points}) in Lemma \ref{28/36}.

(2) Consider the subspace $A_{L_+}'$ of $A_{L_+}$ on which $q_{L_+}$ takes integral values.
Then $(A_{L_+}', q_{L_+})$ is a quadratic form of dimension 7 over $\bbF_2$.  It has a radical
generated by $\tilde{\kappa}/2$.

\medskip

If $\omega\in H^0(X,\Omega_X^2)$, then $\iota^* (\omega_X) =-\omega_X$, and hence
the period domain of the above $K3$ surfaces is given by
\begin{equation}\label{}
\calD = \{ [\omega] \in \bbP(L_-\otimes \bbC) : \la \omega, \omega \ra = 0, \ 
\la \omega, \bar{\omega} \ra > 0 \}
\end{equation}
which is a bounded symmetric domain of type IV and of dimension 12.

\subsection{An isometry of order 4}\label{sigma} 
We shall study the action of $\sigma$ on $H^2(X, {\bbZ})^-$.  Recall that
$$D_4 \cong \{ (x_1,x_2,x_3,x_4) \in {\bbZ}^4 \ \mid \ x_1+x_2+x_3+x_4 \equiv 0 \ (\rm{mod}\ 2) \}.$$
Here we consider the standard inner product on ${\bbZ}^4$ with the {\it negative} sign.
Let $\rho_0$ be the isometry of $D_4$  given by
$$\rho_0(x_1,x_2,x_3,x_4) = (x_2,-x_1,x_4,-x_3).$$
Obviously $\rho_0$ is of order 4 and fixes no non-zero vectors in $D_4$.  
Also an easy calculation shows
that $\rho_0$ acts trivially on $D_4^*/D_4$.  Next let $e,f$ (resp. $e', f'$) be a basis of
$U$ (resp. $U(2)$).   Define the isometry $\rho_1$ of $U \oplus U(2)$ by
$$\rho_1(e) = -e-e', \ \rho_1(f) = f-f',  \ \rho_1(e') = e'+2e, \ \rho_1(f') = 2f-f'.$$
Obviously $\rho_1$ is of order 4, fixes no non-zero vectors in $U \oplus U(2)$ and acts trivially on
the discriminant group of $U \oplus U(2)$.  

Finally let $\rho_2$ be an isometry of $A_1^{\oplus 2}$ with a basis $r_1, r_2$ defined by
$$\rho_2(r_1) = -r_2, \ \rho_2(r_2) = r_1.$$
Now we have an isometry
$\rho = \rho_1 \oplus \rho_0 \oplus \rho_0 \oplus \rho_2$ of order $4$
of $L_{-} = U\oplus U(2) \oplus D_4\oplus D_4\oplus A_1 \oplus A_1$
which fixes no non-zero vectors in $L_-$.   The action of $\rho$ on $L_-^*/L_-$ coincides with the one of $\sigma^*$ on
$H^2(X,\bbZ)_-^*/H^2(X,\bbZ)_-$.

\subsection{A complex ball as the period domain of $K3$ surfaces with an automorphism of order 4}

Now we can define the period domain of the pairs $(X, \sigma)$.  Since $H^0(X,\Omega_X^2) \cong \bbC$,
the period $\omega_X$ is an eigenvector of $\sigma^*$.  First note that
$\rho$ has no non-zero fixed vectors in $L_-\otimes \bbQ$ and hence
the eigenvalues of $\rho$ are $\pm \sqrt{-1}$.  Let
$$V_{\pm} = \{ \omega \in L_-\otimes \bbC : \rho(\omega) =\pm \sqrt{-1}\omega \}$$
both of which have dimension 7. 
Then the period domain of the pairs $(X, \sigma)$ is defined by
\begin{equation}
\calB = \calD \cap \bbP(V_+).
\end{equation}
 If $\omega \in V_{\pm}$, then $\la \omega, \omega \ra =0$.
This shows
$$\calB = \{ [\omega] \in \bbP(V_+) : \la \omega, \bar{\omega} \ra > 0 \}$$ 
and hence $\calB$ is a 6-dimensional complex ball.

Let $r \in L_-$ with $r^2 =\la r, r\ra = -2$.  We define
\begin{equation}
r^{\perp} = \{ \omega \in \calD \ : \ \la \omega , r \ra = 0\}, \ H_r = r^{\perp} \cap \calB, \ \calH = \bigcup_{r\in L_-, r^2=-2} H_r
\end{equation}

\subsection{Proposition}\label{fact2}
{\it There are two types of $(-2)$-vectors $r$ in $L_-$ according to
$r/2 \in L_-^*$ or $r/2 \notin L_-^*$.}

\begin{proof}
See \cite{Kon1}, Lemma 3.3.
\end{proof}

According to the above Proposition, we can write
$$\calH = \calH_n \cup \calH_h$$
where $\calH_n$ (resp. $\calH_h$) is the union of $r^{\perp}$ with $r/2 \notin L_{L_-}^*$
(resp. $r/2 \in L_{L_-}^*$).  We showed that a generic point in $\calH_n$ 
(resp. $\calH_h$) corresponds to a plane quartic with a node 
(resp. a hyperelliptic curve of genus 3) (\cite{Kon1}, Theorems 4.3, 5.3).
We shall study more details of $\calH$ in \ref{heegner1}.

Also we consider the following arithmetic subgroups
\begin{equation}
\Gamma = \{ \gamma \in \O(L_-) : \gamma \circ \rho = \rho \circ \gamma \}, \ 
\tilde{\Gamma} = \Ker \{\Gamma \to \O(q_{L_-}) \}. 
\end{equation}
Then $\Gamma/\tilde{\Gamma}$ is isomorphic to a split extension of $\Sp(6,\bbF_2)$
by $\bbZ/2\bbZ$, where $\bbZ/2\bbZ$ is generated by $\rho$ (e.g. \cite{Kon1}, Lemma 2.2).

\subsection{Proposition}\label{}
{\it The action of $\sigma^*$ on $H^2(X,\bbZ)_-$ is conjugate to $\rho$.}

\begin{proof}
Recall that the action of $\rho$ on $L_-^*/L_-$ coincides with the one of $\sigma^*$ on
$H^2(X,\bbZ)_-^*/H^2(X,\bbZ)_-$.  It follows from Nikulin \cite{N1}, Proposition 1.6.1 that $\rho$ can be extended to an isometry
of $L$ whose action on $L_+$ isomorphic to the one of $\sigma^*$ on $H^2(X,\bbZ)_+$.  
Let $\omega \in \calB\setminus \calH$.  Then the surjectivity of the period map for $K3$ surfaces, there exists a $K3$ surface $Y$
and an isometry $\alpha_Y: H^2(Y,\bbZ)\to L$ with $\alpha_Y(\omega_Y) =\omega$.  Then $\rho$ is represented by an automorphism $\sigma'$ (\cite{Kon1}, Lemma 2.4).  It follows from Nikulin \cite{N2}, Theorem 4.2.2 that the set of fixed points of
$(\sigma')^2$ is a smooth curve $C'$ of genus 3.
We can easily see that the set of fixed points of $\sigma'$ is $C'$, the quotient surface $Y/\la \sigma'\ra$ is a projective 
plane and the branch locus is a smooth plane quartic.  Since the moduli space of plane quartics is connected, the assertion follows.
\end{proof}

\subsection{Proposition}\label{kondo}(\cite{Kon1}, Theorem 2.5)\ 
{\it The moduli space of smooth plane quartics $($ resp. plane quartics with a level $2$ 
structure$)$ is isomorphic to the quotient $(\calB\setminus \calH)/\Gamma$ 
$($ resp. $(\calB\setminus \calH)/(\tilde{\Gamma}\cdot \bbZ/2\bbZ))$.}

\medskip

\section{Hermitian form and reflections}\label{}

\subsection{Hermitian form}\label{hermite}
We consider $L_-$ as a free ${\bbZ}[\sqrt{-1}]$-module $\Lambda$ by
$$(a+b\sqrt{-1})x = ax + b\rho(x).$$
Let 
$$h(x,y) = \sqrt{-1}\langle x, \rho(y) \rangle + \langle x, y \rangle.$$
Then $h(x,y)$ is a hermitian form on ${\bbZ}[\sqrt{-1}]$-module $\Lambda$.
With respect to a ${\bbZ}[\sqrt{-1}]$-basis $(1,-1,0,0)$, $(0,1,-1,0)$ of $D_4$,
the hermitian matrix of $h\mid D_4$ is given by

\begin{equation}\label{}
\begin{pmatrix}-2&1-\sqrt{-1}
\\1+\sqrt{-1}&-2
\end{pmatrix}.
\end{equation}

\noindent
And with respect to a ${\bbZ}[\sqrt{-1}]$-basis $e,e'$ of $U\oplus U(2)$,
the hermitian matrix of $h\mid U\oplus U(2)$ is given by

\begin{equation}\label{}
\begin{pmatrix}0&1+\sqrt{-1}
\\1-\sqrt{-1}&0
\end{pmatrix}.
\end{equation}

\noindent
And with respect to a ${\bbZ}[\sqrt{-1}]$-basis $r$ of $A_1\oplus A_1$,
the hermitian matrix of $h\mid A_1\oplus A_1$ is given by $(-2)$.
Thus the hermitian matrix of $h$ is given by

\begin{equation}\label{}
\begin{pmatrix}0&1+\sqrt{-1}&0&0&0&0&0
\\1-\sqrt{-1}&0&0&0&0&0
\\0&0&-2&1-\sqrt{-1}&0&0&0
\\0&0&1+\sqrt{-1}&-2&0&0&0
\\0&0&0&0&-2&1-\sqrt{-1}&0
\\0&0&0&0&1+\sqrt{-1}&-2&0
\\0&0&0&0&0&0&-2
\end{pmatrix}.
\end{equation}

\noindent
Let 
$$\varphi : \Lambda \to L_-^*$$
be a linear map defined by $\varphi(x) = (x + \rho(x))/2.$
Note that $\varphi((1-\sqrt{-1})x) = \varphi(x - \rho(x)) = x \in L_-$. Hence
$\varphi$ induces an isomorphism
\begin{equation}\label{}
\Lambda/(1-\sqrt{-1})\Lambda \simeq L_-^*/L_-.
\end{equation}

\subsection{Reflections}\label{reflection1}
For $r \in L_-$ with $\langle r, r \rangle = -2$, we define a {\it reflection} 
$$s_r(x) = x + \langle r, x \rangle r$$
which is contained in $\tilde{\O}( L_-) = {\rm Ker}(\O(L_-) \to \O(q_{L_-}))$, but not in $\Aut(\Lambda)\cong \Gamma$.  Also for $\xi \in L_-$ with $\xi^2=-4$ and $\xi/2 \in L_-^*$, we define a reflection $s_{\xi}$ in $\O(L_-)$ by
$$s_{\xi}(x)=x+{1\over 2} \la x, \xi \ra \xi.$$
The reflection $s_{\xi}$ induces a {\it transvection} of $A_{L_-}$ defined by
$$t_{\alpha}(x) = x + 2b_{L_-}(x,\alpha ) \alpha$$
where $\alpha \in A_{L_-}$ is a non-isotropic vector represented by $\xi/2$.

On the other hand, by  considering $r$ as in $\Lambda$, 
we define a {\it reflection} in $\Gamma$ by

\begin{equation}\label{}
R_{r,\epsilon}(x) = x -(1-\epsilon){h(r,x)\over h(r,r)}r
\end{equation}
where $\epsilon \not= 1$ is a 4-th root of unity.  
We can easily see that $R_{r,-1}$ corresponds to the isometry of $L_-$
$$x \to x + \langle r, x \rangle r + \langle \rho(r), x \rangle \rho(r)$$
which coincides with $s_r \circ s_{\rho(r)}$.  Also
$R_{r,\sqrt{-1}}$ corresponds to the isometry of $L_-$
$$x \to x + \langle r, x \rangle (r-\rho(r))/2 + \langle \rho(r), x \rangle (r + \rho(r))/2$$
which induces a transvection $t_{\alpha}$ of $A_{L_-}$ 
where $\alpha \in A_{L_-}$ is a non-isotropic vector $ (r + \rho(r))/2 \ {\rm mod} \ L_-$.
We can easily see that 
$$R_{r,\sqrt{-1}}^2 = s_{r - \rho(r)} \circ s_{r + \rho(r)}.$$
Moreover a direct calculation shows that $R_{r,\sqrt{-1}}^2 = R_{r,-1}$.  We remark that $R_{r,\sqrt{-1}}^2$ acts trivially on
$A_{L_-}$.
It is known that the transvections $t_{\alpha}$ $(\alpha \in A_{L_-}, \ q_{L_-}(\alpha)=1)$ generate $\O(q_{L_-})$.
Thus we have

\subsection{Proposition}\label{}
{\it The natural map 
$$\Gamma \to \O(q_{L_-})$$
is surjective}.
\medskip

\section{Heegner divisors}\label{}

In this section we shall study the discriminant locus $\calH$.

\subsection{Discriminant quadratic form}\label{disc2}

Recall that $L_{+} = A_1(-1) \oplus A_{1}^{\oplus 7}$, $L_{-} = U \oplus U(2) 
\oplus D_4^{\oplus 2} \oplus A_{1}^{\oplus 2}$.  

\medskip
Let $(A_{L_-}, q_{L_-}) = (L_{-}^{*}/L_{-}, q_{L_{-}})$ be the discriminant 
quadratic form of $L_{-}$.  Then $A_{L_-} \simeq (\bbF_2)^{8}$ and
$$q_{L_-} : A_{L_-} \to \bbQ/2\bbZ.$$
Since $L_+$ and $L_-$ are mutually orthogonal complement in the unimodular lattice $L$,
$A_{L_+}$ is canonically isomorphic to $A_{L-}$.
Define
$$A_{L_-}'  = \{ x \in A_{L_-} : q_{L_-}(x) \in \bbZ/2\bbZ \}.$$
Then the quadratic form $q_{L_-} \mid A_{L_-}'$ over $\bbF_2$
has a radical $\la\kappa\ra$ and modulo radical
it defines a symplectic form of dimension 6 over $\bbF_2$ (see Remark \ref{e7-3}, (2)).

\medskip
\noindent
$A_{L_-}$ consists of the following 256 vectors:
\smallskip

Type (00): $ q_{L_-}(x) = 0, \# x = 1, x = 0;$
\smallskip

Type (0): $q_{L_-}(x) = 0, \# x = 63;$
\smallskip	

Type (1): $q_{L_-}(x) = 1, \# x = 63$;
\smallskip

Type (10): $q_{L_-}(x) = 1, \# x = 1.$ $x=\kappa$;
\smallskip

Type (1/2): $q_{L_-}(x) = 1/2, \# x = 56$;
\smallskip

Type (3/2): $q_{L_-}(x) = 3/2, \# x = 72$.
\medskip
\noindent

\medskip
Since $\O(L_-)$ acts transitively on vectors in $A_{L_-}$ with the same type, 
the number of $\tilde{\O}(L_-)$-equivalence classes are 1, 63, 63, 1, 56 or 72 according to
the type (00), (0), (1), (10), (1/2) or (3/2) respectively. 

\subsection{Lemma}\label{e7}
{\it The quadratic space $E_7/2E_7$ is isomorphic to $A_{L_-}'$.}
\begin{proof}
First note that $q_{L_+} = -q_{L_-}$ because $L_-$ is the orthogonal complement of $L_+$ in the unimodular lattice $L$.
On the other hand, the Picard lattice of a del Pezzo surface $S$ of degree 2 is isomorphic to $L_+(1/2)$.  Moreover 
$E_7/2E_7$ is isomorphic to $q_{L_-}\mid A_{L_-}'$ by using the description of $A_{L_+}$ in Lemma \ref{e7-1}.
\end{proof}

\subsection{Remark}\label{e7-5}
Under the isomorphism in Lemma \ref{e7}, 63 discriminant conditions (\ref{discriminant condition}) (resp. 36 conditions 
(\ref{4points})  in Lemma \ref{28/36}) correspond to
63 vectors of type (1) (resp. 72 vectors of type (3/2) modulo the action of $\rho$).
See Remark \ref{e7-3}.

\subsection{Heegner divisors}\label{heegner1}

We now introduce Heegner divisors in $\calD$ and $\calB$.
Let $r \in L_-^*$ with $r^2 < 0$.  We denote by $r^{\perp}$ the orthogonal 
complement in $\calD$: 
$$r^{\perp} = \{ [\omega] \in \calD : \la r, \omega \ra = 0\}.$$
Let $\alpha \in A_{L_-}$ and $n \in \bbQ, \ n < 0$.  Then we difine a {\it Heegner divisor}
$\calD_{\alpha, n}$ in $\calD$ as follows:
$$\calD_{\alpha, n} = \bigcup_{r} r^{\perp}$$
where $r$ moves over the set of all $r \in L_-^*$
satisfying $r\ \mod\ L_- =\alpha$ and $r^2 = n$.

In the later we consider only special cases, that is, $\alpha$ is of type (1) , (10) and $n=-1$,
$\alpha$ is of type (1/2) and $n=-3/2$, and $\alpha$ is of type (3/2) and $n=-1/2$.
In these cases, for simplicity, we denote $\calD_{\alpha, n}$ by $\calD_{\alpha}$:
if $\alpha$ is of type $(1)$ or type $(10)$, $\calD_{\alpha} = \calD_{\alpha, -1}$, if $\alpha$ is of type $(1/2)$, 
$\calD_{\alpha} = \calD_{\alpha, -3/2}$ and if $\alpha$ is of type $(3/2)$, $\calD_{\alpha} = \calD_{\alpha, -1/2}$.
Moreover for $j=1, 10, 1/2$ or $3/2$, we define
$$\calD_j =\bigcup_{\alpha\ {\rm is\ of\ type}\ (j)} \calD_{\alpha}.$$

Let $r\in L_-^*$ with $r^2 < 0$.  Since $\la r, \rho (r)\ra = \la \rho(r), \rho^2(r)\ra =
-\la r, \rho(r)\ra$, $\la r, \rho(r)\ra =0$.  If $\omega \in \calB$, then
$\la r, \omega \ra = \sqrt{-1}\la \rho(r), \omega\ra$.  Thus we have 
$$\la r, \omega \ra = 0 \iff \la \rho(r), \omega \ra= 0 \iff \la r + \rho(r), \omega \ra = 0
\iff \la r - \rho(r), \omega \ra = 0.$$
Therefore
\begin{equation}\label{fact1}
H_r=H_{\rho(r)} = H_{r+ \rho(r)} = H_{r-\rho(r)}
\end{equation}
where $H_r = r^{\perp}\cap \calB$.

\subsection{Lemma}\label{fact3}
{\it Let $r\in L_-$ with $r^2=-2$.  If $r/2\notin L_-^*$, then $(r+\rho(r))/2$ is contained in 
$L_-^*$ and $(r+\rho(r))/2\ {\rm mod}\ L_- $ is of type $(1)$.  If $r/2 \in L_-^*$, then $(r+\rho(r))/2\ {\rm mod} \ L_-$ is of type $(10)$.}
\begin{proof}
The assertion follows from \cite{Kon1}, Lemma 3.3 and its proof.
\end{proof}

\subsection{Lemma}\label{fact4}
{\it Let $r\in L_-$ with $r^2=-2$.}

(1) {\it   The order of the reflection $R_{r,\sqrt{-1}}$ is $4$ and $R_{r,\sqrt{-1}}^2$ is contained in $\tilde{\Gamma}$.}

(2)  {\it If $r/2\notin L_-^*$, then the reflection $R_{r,\sqrt{-1}}$ is not contained in 
$\tilde{\Gamma}\cdot \bbZ/2\bbZ$. }

(3) {\it  If $r/2 \in L_-^*$, then $R_{r,\sqrt{-1}}$ is contained in $\tilde{\Gamma}\cdot \bbZ/2\bbZ$.}
\begin{proof}
A direct calculation shows that $R_{r,\sqrt{-1}}$ is of order 4.  As mentioned in \ref{reflection1}, 
$R_{r,\sqrt{-1}}^2$ acts trivially on $A_{L_-}$.   On the other hand,
$R_{r,\sqrt{-1}}$ induces a transvection $t_{\alpha}$ on $A_{L_-}$ where $\alpha = (r+\rho(r))/2 \ {\rm mod}\ L_-$ (see \ref{reflection1}), and hence it acts non trivially on $A_{L_-}$.
If $r/2\notin L_-^*$ (resp. $r/2\in L_-^*$), then
$\alpha$ is of type (1) (resp. of type (10)) (Lemma \ref{fact3}).  Moreover if $\alpha$ is of type (10), then it acts on $A_{L_-}$ as $\rho$.
Hence the assertion follows.
\end{proof}

For $\alpha \in A_{L_-}$, define a {\it Heegner divisor} in $\calB$ by
$$\calB_{\alpha} = \calB\cap \calD_{\alpha}$$
and denote its image on $\calB/\tilde{\Gamma}\cdot \bbZ/2\bbZ$ by $\calH_{\alpha}$.
It follows from the equation \ref{fact1}, Proposition \ref{fact2}, Lemma \ref{fact3} that
$$\calH = \bigcup_{\alpha} \calB_{\alpha}$$
where $\alpha$ runs over the set of all type $(1), (10), (3/2)$.
Note that if $r\in L_-$ with $r^2=-2$ and $r/2\in L_-^*$ then $(r+\rho(r))/2\ {\rm mod}\ L_- = \kappa$ (Lemma \ref{fact3}).
This implies that $H_r \subset \calB_{\kappa}$.  
Conversely for any $(-1)$-vector $\xi$ in $L_-^*$ with $\xi \ {\rm mod}\ L_- = \kappa$, $\la \xi, \rho(\xi)\ra=0$ and
$\xi +\rho(\xi)$ is a $(-2)$-vector in $L_-$.  Obviously $(\xi +\rho(\xi))/2$ is of type $(3/2)$.
Thus we have
$$\calB_{\kappa} = \calB \cap \calD_{10} = \bigcup_{\alpha} \calB_{\alpha}$$
where $\alpha$ runs over the set of vectors of type (3/2).  
We now conclude:

\subsection{Lemma}\label{}
{\it $\calH/(\tilde{\Gamma}\cdot \bbZ/2\bbZ)$ consists of $63$ components $\calH_{\alpha} \ (\alpha$ is of type $(1))$ and
$36$ components $\calH_{\beta}$ $(\beta$ is of type $(3/2))$. }

\subsection{Lemma}\label{branch}
{\it  For any $\alpha$, $\calH_{\alpha}$ is a branch divisor  of the covering $\calB \to \calB/(\tilde{\Gamma}\cdot \bbZ/2\bbZ)$. 
 If $\alpha$ is of type $(1)$, then the branch degree is two.  If $\alpha$ is of type $(3/2)$, the branch degree is four.}
\begin{proof}
We use Lemma \ref{fact4}.
 If $\alpha$ is of type $(1)$, then there exists $r \in L_-$ with $r^2=-2$ such that  $\alpha = (r+\rho(r))/2\ {\rm mod}\ L_-$.  
Then $R_{r,\sqrt{-1}}^2$ is contained in $\tilde{\Gamma}$ and fixes $H_r$.
On the other hand, $\alpha$ is type $(3/2)$,  there exists $r \in L_-$ with $r^2=-2$ such that  $\alpha = r/2\ {\rm mod}\ L_-$.
Then $R_{r,\sqrt{-1}}$ is contained in $\tilde{\Gamma}\cdot \bbZ/2\bbZ$ and fixes $H_r$.  Thus 
the assertion follows.
\end{proof}

\medskip

\subsection{Generalized del Pezzo surfaces}\label{gendel}
Let  $\pi : \hat{P}_2^7 \to P_2^7$ be the blow up of $P_2^7$ along the 21 subvarieties of codimension 2 corresponding to
21 roots $e_i-e_j$ (see Lemma \ref{blowup}).  
Let $(\hat{P}_2^7)_0$ be the open set of $\hat{P}_2^7$ consisting of the inverse image of ordered 7 points in general position.
Let $(\hat{P}_2^7)_1$ be the complement of the inverse images of 36 subvarieties corresponding to 
the conditions given in Lemma \ref{28/36}.
Recall that $\Sp(6,\bbF_2)$ naturally acts on $P_2^7$ as birational automorphisms (\cite{DO}, Chap. VII, \S 4).
It follows from Proposition \ref{kondo} that  $(\hat{P}_2^7)_0$ is isomorphic to
$(\calB\setminus \calH)/(\tilde{\Gamma}\cdot \bbZ/2\bbZ)$.  
We denote this isomorphism by $p$.  Note that $p$ is $\Sp(6,\bbF_2)$-equivariant. 

Let $p_1,..., p_7$ be an ordered 7 points in $\bbP^2$ representing a point in $(\hat{P}_2^7)_1\setminus (\hat{P}_2^7)_0$.
For simplicity, we consider generic case and we may assume that $p_1,..., p_6$ lies on a conic $Q$.   Then after blowing up at 7 points,
we get a {\it generalized} del Pezzo surface $S$ containing a $(-2)$-curves, i.e., the proper transform of $Q$.  
The anti-canonical map of $S$ gives a double cover of
$\bbP^2$ branched along a nodal quartic curve $C$.  In \cite{Kon1}, \S 4, we studied the $K3$ surface $X$ which is the 
minimal model of the 4-cyclic cover of $\bbP^2$ branched along $C$.  The Picard lattice of a generic $X$ is isomorphic to
$U\oplus A_1^{\oplus 8}$.  Moreover $X$ has an elliptic fibration with 8 singular fibers of type $III$ in the sense of Kodaira.
This elliptic fibration is induced from the pencil of lines in $\bbP^2$ through the node of $C$.  8 singular fibers correspond to
two tangent lines of $C$ at the node and six tangent lines of $C$ at some smooth points (note that the normalization of $C$ is a
hyperelliptic curve of genus 2).  On the other hand, this elliptic fibration is also induced from the pencil of lines on $\bbP^2$
through $p_7$.  8 singular fibers correspond to six lines through $p_1,...,p_6$ and two tangent lines $l_1, l_2$ of the conic $Q$.
To get $X$ from $S$, we need a blow up at two points $q_1, q_2$ at which $l_1,l_2$ tangent to $Q$.
Note that  the $K3$ surface $X$ has new algebraic cycles.  Two exceptional curves over $q_1, q_2$ gives
a sublattice $A_1\oplus A_1$ in $L_-$.  
Denote by $r_1, r_2$ a generator of 
$A_1\oplus A_1$.  Then $r_2 = -\rho(r_1)$ and $(-1)$-vector $(r_1+r_2)/2$ in $L_-^*$ corresponds to 
$(2\tilde{e}_0 - \tilde{e}_1 - \cdot \cdot \cdot - \tilde{e}_6)/2$ in Lemma \ref{e7-1} under the isomorphism 
between $L_+^*/L_+ \cong L_-^*/L_-$.  The sum of these two $(-1)$-vectors is represented by the class of $(-2)$-curve corresponding to
the conic $Q$.
Since $\omega_X$ is perpendicular to $A_1\oplus A_1$, 
the period of the $K3$ surface $X$ is contained in $\calB_{\alpha}$ where $\alpha = (r_1+r_2)/2 \ {\rm mod} \ L_-$.
Now we conclude:

\subsection{Proposition}\label{kondo1} 
{\it $p$ can be extended a holomorphic map 
$$\hat{p} : (\hat{P}_2^7)_1 \to (\calB\setminus \calH_h)/(\tilde{\Gamma}\cdot \bbZ/2\bbZ)$$ 
which is $\Sp(6,\bbF_2)$-equivariant and isomorphic on the complement of subvarieties of codimension $2$.}

\begin{proof}
Let $P$ be a component from 63 divisors in $\hat{P}_2^7$.  Take a general point of $P$ and consider a corresponding
7 points on $\bbP^2$.  As mentioned above, we have a $K3$ surface associated to these 7 points whose period is contained
in $\calB_{\alpha}$.  Thus,
by using the theory of simultaneous resolutions of singularities, $p$ can be extended to a holomorphic map $\hat{p}$ to a general
point of $P$.
It follows from the proof of Theorem 4.3 in \cite{Kon1} that $\hat{p}$ is injective on a general point of $P$.
Since $p$ is equivariant under the action of 
$\Sp(6,\bbF_2)$, its extension $\hat{p}$ is $\Sp(6,\bbF_2)$-equivariant, too.
\end{proof}

\medskip

\section{Weil representtion}\label{}

In this section we study an action of $\SL(2, \bbZ)$ on the group ring $\bbC[A_{L_-}]$
called Weil representation.   The following Lemmas \ref{char}, \ref{key} will be used in the next section to show the existence of
a $15$-dimensional linear system of automorphic forms on $\calB$.
For simplicity, we sometimes denote the discriminant quadratic form $q_{L_-}$ by $q$
and the discriminant bilinear form $b_{L_-}$ by $ b$.
In the following we denote by $S, T$ a generator of $\SL(2, \bbZ)$:

\begin{equation}\label{}
T =
\begin{pmatrix}1&1
\\0&1
\end{pmatrix},\quad
S =
\begin{pmatrix}0&-1
\\1&0
\end{pmatrix}.
\end{equation}

Let $\rho$ be the Weil representation of $\SL(2, \bbZ)$ on $\bbC[A_{L_-}]$ defined by:

\begin{equation}\label{Weil}
\rho(T)(e_{\alpha}) = e^{\pi\sqrt{-1}\ q(\alpha)} e_{\alpha}, \quad
\rho(S)(e_{\alpha}) = {\sqrt{-1} \over 16} \sum_{\delta} 
e^{-2\pi\sqrt{-1}\ b( \delta, \alpha )} e_{\delta}.
\end{equation}

The action of $\SL(2,\bbZ)$ on $\bbC[A_{L_-}]$ factorizes to the one of $\SL(2,\bbZ/4\bbZ)$.
The conjugacy classes of $\SL(2, \bbZ/4\bbZ)$ consist of
$\pm E, \pm S, \pm T, \pm T^{2}, ST, (ST)^2$.  
Let $\chi_{i}$ $(1 \leq i \leq 10)$ be the characters of irreducible
representations of $\SL(2, \bbZ/4\bbZ)$.  For the convenience of the reader we give the
character table of $\SL(2, \bbZ/4\bbZ)$ in the Table 1.

\begin{table}[h]
\[
\begin{array}{rlllllllllllllllllllllll}
 & E&-E&S&-S&T&-T&T^2&-T^2&ST&(ST)^2\\
\chi_1&1&1&1&1&1&1&1&1&1&1\\
\chi_2&1&1&-1&-1&-1&-1&1&1&1&1\\
\chi_3&1&-1&\sqrt{-1}&-\sqrt{-1}&\sqrt{-1}&-\sqrt{-1}&-1&1&-1&1\\
\chi_4&1&-1&-\sqrt{-1}&\sqrt{-1}&-\sqrt{-1}&\sqrt{-1}&-1&1&-1&1\\
\chi_5&2&2&0&0&0&0&2&2&-1&-1\\
\chi_6&2&-2&0&0&0&0&-2&2&1&-1\\
\chi_7&3&3&1&1&-1&-1&-1&-1&0&0\\
\chi_8&3&3&-1&-1&1&1&-1&-1&0&0\\
\chi_9&3&-3&-\sqrt{-1}&\sqrt{-1}&\sqrt{-1}&-\sqrt{-1}&1&-1&0&0\\
\chi_{10}&3&-3&\sqrt{-1}&-\sqrt{-1}&-\sqrt{-1}&\sqrt{-1}&1&-1&0&0\\
\end{array}
\]
\caption{}
\end{table}

On the other hand, 
for each $u \in A_{L_-}$, we denote by $m_{0}$ (resp. $m_{1}$) the number of vectors 
$v \in A_{L_-}$ with $b_{L_-}(u, v) \equiv 0$ 
(resp. $1/2$).   Then $m_{0}, m_{1}$ are given in the Table 2.

\begin{table}[h]
\[
\begin{array}{rlllllllllllllllllllllll}
u& 00&00&00&00&00&00&0&0& 0&0&0&0&1&1&1&1&1&1\\
v&00&0&1&10&1/2&3/2&00&0&1&10&1/2&3/2&00&0&1&10&1/2&3/2\\
m_0&1&63&63&1&56&72&1&31&31&1&24&40&1&31&31&1&32&32\\
m_1&0&0&0&0&0&0&0&32&32&0&32&32&0&32&32&0&24&40\\
u& 10&10&10&10&10&10&1/2&1/2& 1/2&1/2&1/2&1/2&3/2&3/2&3/2&3/2&3/2&3/2\\
v&00&0&1&10&1/2&3/2&00&0&1&10&1/2&3/2&00&0&1&10&1/2&3/2\\
m_0&1&63&63&1&0&0&1&27&36&0&28&36&1&35&28&0&28&36\\
m_1&0&0&0&0&56&72&0&36&27&1&28&36&0&28&35&1&28&36\\
\end{array}
\]
\caption{}
\end{table}

\subsection{Lemma}\label{char}
{\it Let $\chi$ be the character of the representation of
$\SL(2, \bbZ/4\bbZ)$ on $\bbC[A_{L_-}]$.
Let $\chi = \sum_{i} m_{i} \chi_{i}$ be
the decomposition into irreducible characters.  Then
$$\chi = 7\chi_3 + 15\chi_4 + 21\chi_6 + 28\chi_9 + 36\chi_{10}.$$}

\begin{proof}
By using the Table 2, we can see that
the traces of conjugacy classes of $\SL(2, \bbZ/4\bbZ)$ are as follows:
$$tr(E) = 2^8, \ tr(-E) = -2^8, \ tr(S) = 0, \ tr(-S) =0, \ tr(T) = -16\sqrt{-1},$$
$$tr(-T) = 16\sqrt{-1}, \ tr(T^{2}) = 0, \ tr(-T^2) = 0, \ tr(ST) = -1, \ tr((ST)^2) = 1.$$
The assertion now follows from the Table 1.
\end{proof}

\subsection{Lemma}\label{totally}
(i)   {\it The number of totally isotropic subspaces of 
dimension $2$ in $A_{L_-}$ is $315$}.

(ii) {\it The number of totally isotropic subspaces of 
dimension $3$ in $A_{L_-}$ is $135$.}

(iii) {\it Let $e$ be a non-zero isotropic vector.  Then there exist $15$ totally isotropic subspaces of dimension $2$ and $35$ totally isotropic subspaces of dimension $3$ containing $e$.
There exist $3$ totally isotropic subspaces of dimension $3$ containing a fixed totally isotropic subspaces of dimension $2$.}
\begin{proof}
The number of totally isotropic subspaces of dimension $2$ or dimension $3$ is given by
$${(2^{6} - 1)2(2^{4} - 1)\over (2^{2} - 1)(2^{2} - 2)} = 315$$
or
$${(2^{6} - 1)2(2^{4} - 1)2^2(2^2-1)\over (2^{3} - 1)(2^{3} - 2)(2^3-2^2)} = 135$$
respectively.  The assertion (iii) now follows from the facts that the number of non-zero isotropic vectors is 63 and
each totally isotropic subspace of dimension 3 contains exactly 7 totally isotropic subspaces of dimension 2.
\end{proof}

\smallskip
\subsection{Definition-Remark}\label{rep}
It follows from Lemma \ref{char} that there exists a 15-dimensional subspace of 
$\bbC[A_{L_-}]$ on which $\SL(2, \bbZ/4\bbZ)$ acts with character $\chi_4$.
We denote by $W$ this 15-dimensional representation.

Note that $\O(q_{L_-})$ naturally acts on $\bbC[A_{L_-}]$.
Since the action of $\SL(2, \bbZ/4\bbZ)$ on $\bbC[A_{L_-}]$ and that of $\O(q_{L_-})$ commute,
$\O(q_{L_-})$ acts on $W$.  One can prove that this 15-dimensional representation of $\O(q_{L_-})$ is irreducible (e.g. see Dolgachev, Ortland \cite{DO}, Chapter IX, the proof of Proposition 9).

\subsection{Definition}\label{invariant}
Let $I$ be a maximal totally isotropic
subspace  of $A_{L_-}$ and let $V=\la I, \kappa\ra$ be a subspace of $A_{L_-}$ generated by $I$ and $\kappa$.  
Let $\alpha \in A_{L_-}$ with 
$q_{L_-}(\alpha) = 3/2$, $b_{L_-}(\alpha, c) = 0$ for any $c \in I$ 
($\alpha$ is unique modulo $V$).
Define
$$M_{+} = \{ \alpha + c : c \in I\}, \quad M_{-} = \{ \alpha + c + \kappa : c \in I\},$$
and
$$\theta_V = \sum_{\beta \in M_{+}} e_{\beta} - \sum_{\beta \in M_{-}}
e_{\beta} \in \bbC[A_{L_-}].$$

\subsection{Lemma}\label{key}
(i)  \ $\rho(S)(\theta_{V}) = -\sqrt{-1} \theta_{V}$, \quad  $\rho(T)(\theta_{V}) = -\sqrt{-1} \theta_{V}.$
{\it In particular $\theta_V$ is contained in $W$.}
\smallskip

(ii) \ {\it  For $a \in V$ with $q_{L_-}(a) = 1$, 
$t_{a}(\theta_{V}) = -\theta_{V}$ where $t_{a}$ is the transvection associated 
to $a$.}

\begin{proof}
(i)  If $\beta \in M_{\pm}$, then $q_{L_-}(\beta)=-1/2$, and hence 
$\rho(T)(\theta_{V}) = -\sqrt{-1} \theta_{V}.$  Next by definition,
$$\rho(S)(\theta_{V}) = {\sqrt{-1} \over 16} \sum_{\delta} (\sum_{\beta\in M_+}
e^{-2\pi\sqrt{-1}\ b(\delta, \beta)} - \sum_{\beta\in M_-} e^{-2\pi\sqrt{-1}\ 
b(\delta, \beta)}) e_{\delta}.$$
We denote by $c_{\delta}$ the coefficient of $e_{\delta}$.  If $\delta \in M_+$, then
$b_{L_-}(\delta, \beta) = -1/2$ for $\beta \in M_+$ and 
$b_{L_-}(\delta, \beta) \in \bbZ$ for $\beta \in M_-$.  
Hence $c_{\delta} = -2^3-2^3=-2^4$.
Similary if $\delta \in M_-$, then $c_{\delta} = 2^4$.

Now assume $\delta \notin M_{\pm}$.  If $\delta \in V$, we can easily see that $c_{\delta} =0$.
Hence we assume $\delta \notin V$.
First consider the case 
$b_{L_-}(\delta, \kappa) \in \bbZ$.
Since $V^{\perp}=V$, there exists $\gamma \in V$ such that 
$b_{L_-}(\gamma, \delta) \notin \bbZ$.  In this case 
$I = \delta^{\perp} \cap I \cup \{ \gamma + a : a \in \delta^{\perp} \cap I \}$.  
This implies that
$$\sum_{\beta\in M_+}
e^{-2\pi\sqrt{-1} \ b(\delta, \beta)} = \sum_{\beta\in M_-} e^{- 2\pi\sqrt{-1}\ b(\delta, \beta)}=0.$$
Finally if $b_{L_-}(\delta, \kappa) \notin \bbZ$, then $\delta = \alpha + \delta'$,
$b_{L_-}(\delta', \kappa) \in \bbZ$.  Then this case reduces to the previous case.

(ii) The transvection $t_a$ interchanges $M_+$ and $M_-$ and hence the assertion follows.
\end{proof}

\subsection{Lemma}\label{linear-rel}
{\it Let $A$ be a totally isotropic subspace of dimension $2$ and let $I_i$ $(i=1,2,3)$ be
totally isotropic subspaces of dimension $3$ contining $A$.  Let $V_i = \la I_i, \kappa\ra$.
Then 
$$\theta_{V_1} \pm \theta_{V_2} \pm \theta_{V_3} = 0$$
where the sign $\pm$ are taken for suitable one.}

\begin{proof}
We fix a decomposition
$$A_{L_-} = u_1\oplus u_2\oplus u_3\oplus q_1 \oplus q_2$$
where $u_i$ is a hyperbolic plane and both $q_1$ and $q_2$ are the discriminant  quadratic form of $A_1$.  
Let $e_i, f_i$ be a basis of $u_i$ $(i=1,2)$ with $q_{L_-}(e_i)=q_{L_-}(f_i)=0$, $b_{L_-}(e_i, f_i)=1/2$
and $\alpha_i$ is a generator of $q_i$ with $q_{L_-}(\alpha_i)=3/2$.  Then $\kappa = \alpha_1 +\alpha_2$ and
$$A'_{L_-} = u_1\oplus u_2\oplus u_3\oplus \la \kappa\ra.$$
We may assume that $A$ is generated by $e_1, e_2$.  Then
$$I_1 = \la e_1, e_2, e_3\ra, I_2=\la e_1,e_2,f_3\ra, I_3=\la e_1,e_2,e_3+f_3+\kappa\ra.$$
We take $\alpha_1$ as $\alpha$ for $I_1, I_2$ and $\alpha_1+e_3$ for $I_3$.
An easy calculation shows that $\theta_{V_1} - \theta_{V_2} = \theta_{V_3}$. 
\end{proof}

\subsection{Lemma}\label{cubic}
 {\it Let $e$ be a non-zero isotropic vector.  Among $15$ totally isotropic subspaces of dimension $3$ containing $e$, there are six $I_j$ $(j=1,...,6)$ such that for other $I$, $\theta_{V}$ $(V=\la I, \kappa\ra)$ is a linear combination of $\theta_{V_j}$ $(V_j=\la I_j, \kappa\ra)$, and
the set of isotropic vectors in $I_1 \cup I_2\cup I_3$ coincides with 
that of $I_4 \cup I_5\cup I_6$ including multiplicities.}

\begin{proof}
We use the same notation as in Lemma \ref{linear-rel} and assume that
$e = e_1$.  Then by using Lemma \ref{linear-rel}, we can see that the following $I_1,..., I_6$ satisfies that
$\theta_V$ is a linear combination of $\theta_{V_1},..., \theta_{V_6}$:
$$I_1=\la e_1, e_2, e_3 \ra, \ I_2=\la e_1, f_2, f_3\ra, \ I_3=\la e_1, e_2+f_3, f_2+e_3\ra,$$
$$I_4 =\la e_1,e_2,f_3\ra, \ I_5=\la e_1, f_2, e_3\ra, \ I_6=\la e_1, e_2+e_3, f_2+f_3\ra.$$
We can easily see that both $I_1 \cup I_2\cup I_3$ and $I_4 \cup I_5\cup I_6$
contain $0$ and $e_1$ with multiplicity 3 and other isotorpic vectors with multiplicity 1.

\end{proof}

\medskip

\section{Automorphic forms: additive liftings}\label{}

In this section we shall construct a linear system of meromorphic automorphic forms by
applying Borcherds theory \cite{B1}.

Let $\rho : \SL(2,\bbZ) \to \GL(\bbC[A_{L_-}])$ be the Weil representation (\ref{Weil}).
A holomorphic map 
$$f : H^+ \to \bbC[A_{L_-}]$$
 is called a {\it vector valued modular form of weight $k$
with respect to }$\rho$ if
$$f(M\tau)=\rho(M)(c\tau+d)^kf(\tau)$$
for any $M = \begin{pmatrix}a&b
\\c&d
\end{pmatrix}
\in \SL(2,\bbZ)$.  For any $\theta \in W$, define
$$f_{\theta}(\tau) = \eta(\tau)^{-6}\theta$$
where $\eta$ is the Dedekind eta function. 
Then it follows that 
$$f_{\theta}(\tau +1) = \eta(\tau +1)^{-6}\theta = -\sqrt{-1}\eta(\tau)^{-6}\theta =\rho(T)f_{\theta}(\tau),$$
$$f_{\theta}(-1/\tau) = \eta(-1/\tau)^{-6}\theta = -\sqrt{-1}\tau^{-3}
\eta(\tau)^{-6}\theta =\tau^{-3}\rho(S)f_{\theta}(\tau),$$
that is, $f_{\theta}(\tau)$ is a vector valued modular form of weight $-3$ with respect to 
$\rho$.  Thus we have a $15$-dimensional space $\tilde{W}$
of modular forms of weight $-3$ and of type $\rho$.

By applying Borcherds \cite{B1}, Theorem 14.3, we have a $\O(q_{L_-})$-equivariant map
$$F: W \to \tilde{W} \to W_{2}(\calD,\tilde{\O}(L_-))$$
where $W_{2}(\calD,\tilde{\O}(L_-))$ is the space of meromorphic automorphic forms of weight 2 on 
a 12-dimensional bounded symmetric domain $\calD$
associated to $L_-$ with respect to $\tilde{\O}(L_-)= \Ker\{\O(L_-)\to \O(q_{L_-})\}$.
For $\theta_V$ in Definition \ref{invariant}, 
we denote $f_{\theta_V}$ by $f_V$ for simplicity, and 
define $F_V= F(\theta_V)$.  Let $\alpha \in V$ be a non isotropic vector.  Take $r \in L_-$ with $r^2=-2$ and 
$\alpha = (r + \rho(r))/2\ {\rm mod}\ L_-$ .  Then the reflection $s_{(r+\rho(r))}$ fixes the hyperplane $(r+\rho(r))^{\perp}$ 
and it induces a transvection $t_{\alpha}$ on $A_{L_-}$ (see \ref{reflection1}).
It now follows from Lemma \ref{key}, (ii) and the $\O(q_{L_-})$-equivariantness of $F$ that
$F_V$ vanishes along 8 Heegner divisors associated  to non isotoropic vectors in $V$.  
Note that the negative power of the Fourier expansion of $f_V$ is $q^{-1/4}\theta_V$
($q=e^{2\pi\sqrt{-1}\tau}$).
It follows from Borcherds \cite{B1}, Theorem 14.3 (also see its proof in page 555) 
that $F_V$ has poles with multiplicity 2 along Heegner divisors associated to 
$16$ vectors with norm $-1/2$ appeared in $\theta_V$.

\subsection{Lemma}\label{}
{\it The additive lifting $F$ is injective.}

\begin{proof}
The map $F$ is $\O(q_{L_-})$-equivariant and $W$ is an irreducible representation of $\O(q_{L_-})$.
Since $F_V$ has a pole, $F_V\not= 0$.  The assertion now follows from the Schur's Lemma.
\end{proof}

We now conclude:

\subsection{Theorem}\label{main1}
{\it For each $V$, $F_V$ is a meromorphic  automorphic form on $\calD$ of weight $2$
vanishing at least along Heegner divisors associated to $8$ non isotropic vectors in $V$
and with poles of multiplicity $2$ along Heegner divisors associated to $16$ vectors with norm $3/2$.
The product of $135$ $F_V$'s is a meromorphic automorphic form on $\calD$ of weight $270$ vanishing along the Heegner divisor $\calD_1$ with at least multiplicity $15$ and along the Heegner divisor $\calD_{10}$ with at least multiplicity $135$ 
and having poles along the Heegner divisor $\calD_{3/2}$ with multiplicity $60$ exactly. }

\begin{proof}
Each $V$ contains 7 non-isotropic vectors of type (1)  and the number of vectors of type (1) is 63.  Hence
the product vanishes along each Heegner divisor of type (1) with at least multiplicity
$${135\times 7 \over 63} = 15.$$
Each $V$ contains $\kappa$ and hence the product vanishes along the Heegner divisor $\calD_{10}$ with at least multiplicity $135$.
Finally the pole order of the product along $\calD_{3/2}$ is given by
$${135\times 16 \times 2 \over 72} = 60.$$
\end{proof}

\medskip

\section{Automorphic forms: Borcherds products}\label{}

In the previous section we construct a linear system of automorphic forms.
To determind the divisor of $F_V$ we need an another automorphic form with
known zeros and poles.  In this section we construct a suitable such automorphic form
by applying Borcherds theory \cite{B1}, \cite{B2}, \cite{F}.

\subsection{Obstruction}\label{}

Borcherds products are meromorphic automorphic forms on $\calD$ 
whose zeros and poles lie on Heegner divisors.  To show the existence of some 
Borcherds products, we introduce the {\it obstraction space}
consisting of all vector valued elliptic modular forms $\{ f_{\alpha} \}_{\alpha \in A_{L_-}}$ of weight $(2+12)/2=7$ and with respect to
the dual representation of $\rho$:

\begin{equation}\label{}
f_{\alpha}(\tau + 1) = e^{-\pi \sqrt{-1}\ q(\alpha)}f_{\alpha}(\tau), \quad
f_{\alpha}(-1/\tau) = -{\sqrt{-1} \tau^{7} \over 2^{4}} 
\sum_{\beta} e^{2\pi \sqrt{-1} \ b(\alpha, \beta)} f_{\beta}(\tau).
\end{equation}

We shall apply the next theorem to show the existence of 
some Borcherds products (see Theorem \ref{multiplicative}, Corollary \ref{multi}).

\subsection{Theorem}\label{freitag}(Borcherds \cite{B2}, Freitag \cite{F}, Theorem 5.2)
{\it A linear combination
$$\sum_{\alpha \in A_{L_-}, n<0} c_{\alpha , n} \calD_{\alpha, n}, \ c_{\alpha, n} \in \bbZ$$
is the divisor of a meromorphic automorphic form on $\calD$ 
of weight $k$ if for every cusp form 
$$f = \{f_{\alpha}(\tau)\}_{\alpha\in A_{L_-}}, 
\ f_{\alpha}(\tau) = \sum_{n \in \bbQ} a_{\alpha, n} e^{2\pi \sqrt{-1} n \tau}$$ 
in the obstruction space, the relation
$$\sum_{\alpha \in A_{L_-}, n<0} a_{\alpha, -n/2}c_{\alpha, n} = 0$$
holds.  In this case the weight $k$ is given by
$$k = \sum_{ \alpha \in A_{L_-}, n\in \bbZ} b_{\alpha, n/2}c_{\alpha, -n}$$
where $b_{\alpha, n}$ are the Fourier coefficients of the Eisenstein series in the obstruction 
space with the constant term
$b_{0, 0} = -1/2$ and $b_{\alpha, 0} = 0$ for $\alpha \not= 0$.}

\medskip

In the following we shall study the divisors $\sum_{\alpha \in A_{L_-}, n<0} c_{\alpha , n} \calD_{\alpha, n}$
where $c_{\alpha, n}$ depends only on the type of $\alpha$.
We denote by 
$$h_{00},\ h_{0},\ h_{1},\ h_{10},\ h_{1/2},\ h_{3/2}$$
the sum of $f_{\alpha}$'s according to their types.  Then we need to consider the 6-dimensional representation $V$  given by:

\begin{equation}\label{}
\rho^*(T) =
\begin{pmatrix}1&0&0&0&0&0
\\0&1&0&0&0&0
\\0&0&-1&0&0&0
\\0&0&0&-1&0&0
\\0&0&0&0&-\sqrt{-1}&0
\\0&0&0&0&0&\sqrt{-1}
\end{pmatrix};
\end{equation}

\begin{equation}\label{}
\rho^*(S) = - {\sqrt{-1} \over 16}
\begin{pmatrix}1&63&63&1&56&72
\\1&-1&-1&1&-8&8
\\1&-1&-1&1&8&-8
\\1&63&63&1&-56&-72
\\1&-9&9&-1&0&0
\\1&7&-7&-1&0&0
\end{pmatrix}.
\end{equation}

\subsection{Lemma}\label{dim}
{\it The dimension of the space of modular forms of weight $7 $ and of 
type $\rho^{*}$ is $4$. The dimension of the space of Eisenstein forms of weight $7$ and of type $\rho^{*}$ is $2$. }

\begin{proof}
In general, the dimension of the space of modular forms of weight $k > 2$ and of type $\rho^*$
is given by
$$d + dk/12 - \alpha (e^{k\pi\sqrt{-1}/2} \rho^*(S)) - \alpha ((e^{k\pi\sqrt{-1}/3}\rho^*(ST))^{-1}) - \alpha (\rho^*(T))$$
(\cite{B2}, section 4, \cite{F}, Proposition 2.1).
Here $d=\dim \{ x \in V : \rho^*(-E)x = (-1)^kx\}$ and 
$$\alpha (A)=\sum_{\lambda} \alpha$$ 
where $\lambda$ runs through all eigenvalues of $A$ and $\lambda = e^{2\pi\sqrt{-1}\alpha}$, $0\leq \alpha <1.$

In our situation, $k=7$ and $d= \dim (V) = 6$.  An elementary calculation shows that
$$\alpha (e^{7\pi\sqrt{-1}/2} \rho^*(S)) = 3/2, \ \alpha ((e^{k\pi\sqrt{-1}/3}\rho^*(ST))^{-1}) = 2, \ \alpha (\rho^*(T)) = 2.$$
On the other hand, the space of Eisenstein series is isomorphic to the subspace of $V$ given by
$$\rho^*(T)(x) =x, \ \rho^*(-E)(x) = (-1)^kx$$
(see Remark 2.2 in \cite{F}).
Thus we have the assertion. 
\end{proof}

\subsection{Eisenstein forms}\label{eisenstein}
In the following we shall calculate Eisenstein forms $\{ h_{\alpha}\}$ 
of weight 7 and of type $\rho^*$.
Then $\{ h_{\alpha}\}$ should satisfy the following:

\begin{equation}\label{}
\begin{cases}
h_{00}(\tau + 1) = h_{00}(\tau), \ h_{0}(\tau + 1) = h_{0}(\tau), 
\ h_{1}(\tau + 1) = -h_{1}(\tau),\cr 
h_{10}(\tau + 1) = -h_{10}(\tau),\ h_{1/2}(\tau + 1) = -\sqrt{-1} h_{1/2}(\tau),\ 
h_{3/2}(\tau + 1) = \sqrt{-1} h_{3/2}(\tau),\cr
h_{00}(-1/\tau) = {-\sqrt{-1} \tau^{7} \over 16} (h_{00} + h_{0} + h_{1} +
h_{10} + h_{1/2} + h_{3/2}),\cr
h_{0}(-1/\tau) = {-\sqrt{-1} \tau^{7} \over 16} (63 h_{00} - h_{0} - h_{1} +
63 h_{10} - 9 h_{1/2} + 7 h_{3/2}),\cr
h_{1}(-1/\tau) = {-\sqrt{-1} \tau^{7} \over 16} (63 h_{00} - h_{0} - h_{1} +
63 h_{10} + 9 h_{1/2} - 7 h_{3/2}),\cr
h_{10}(-1/\tau) = {-\sqrt{-1} \tau^{7} \over 16} (h_{00} + h_{0} + h_{1} +
h_{10} - h_{1/2} - h_{3/2}),\cr
h_{1/2}(-1/\tau) = {-\sqrt{-1} \tau^{7} \over 16} (56 h_{00} - 8 h_{0} + 
8h_{1} - 56 h_{10}),\cr
h_{3/2}(-1/\tau) = {-\sqrt{-1} \tau^{7} \over 16} (72 h_{00} + 8 h_{0} - 
8 h_{1} - 72 h_{10}).\cr
\end{cases}
\end{equation}

Let 
$$E_{1} = G_{7}^{(0,1)}(\tau),\ E_{2} = G_{7}^{(1,0)}(\tau),\ E_{3} = G_{7}^{(1,1)}(\tau),$$
$$E_{4} = G_{7}^{(1,2)}(\tau),\ E_{5} = G_{7}^{(1,3)}(\tau),\ E_{6} = G_{7}^{(2,1)}(\tau)$$
be the Eisenstein series of weight 7 and level 4.  Here $(0,1), (1,0), (1,1), (1,2), (1,3), (2,1) \in (\bbZ/2\bbZ)^2$.
Then their Fourier expansion are given as
follows (Koblitz, \cite{Kob}, Chap. III, \S 3, Proposition 22):

\begin{equation}\label{}
\begin{cases}
E_{1} = {(2\pi)^{7} \over 2\cdot 7!}B_{7}(1/4) + {(-2\pi \sqrt{-1})^{7} 2\sqrt{-1}
\over 4^{7}\cdot 6!} q + \cdot \cdot \cdot;\cr
E_{2} = {(-2\pi \sqrt{-1})^{7} \over 4^{7}\cdot 6!}(q^{1/4} + 2^{6} q^{1/2} +
(3^{6} - 1) q^{3/4} + 4^{6} q + \cdot \cdot \cdot) ;\cr
E_{3} = {(-2\pi \sqrt{-1})^{7} \over 4^{7}\cdot 6!}(\sqrt{-1} q^{1/4} - 2^{6} q^{1/2} 
- \sqrt{-1}(3^{6} - 1) q^{3/4} + 4^{6} q + \cdot \cdot \cdot) ;\cr
E_{4} = {(-2\pi \sqrt{-1})^{7} \over 4^{7}\cdot 6!}(- q^{1/4} + 2^{6} q^{1/2} -
(3^{6} - 1) q^{3/4} + 4^{6} q + \cdot \cdot \cdot) ;\cr
E_{5} = {(-2\pi \sqrt{-1})^{7} \over 4^{7}\cdot 6!}(-\sqrt{-1} q^{1/4} - 2^{6} q^{1/2} 
+ \sqrt{-1}(3^{6} - 1) q^{3/4} + 4^{6} q + \cdot \cdot \cdot) ;\cr
E_{6} = {(-2\pi \sqrt{-1})^{7} \over 4^{7}\cdot 6!}(2\sqrt{-1} q^{1/2} 
+ 0 \cdot q + \cdot \cdot \cdot)\cr
\end{cases}
\end{equation}

where $B_7$ is the Bernoulli polynomial.
The action of $S, T$ on $E_i$ is as follows:

$$T(E_1)=E_1,\ T: E_2\to E_3\to E_4\to E_5 \to E_2,\ T(E_6)=-E_6;$$
$$S: E_1\to E_2\to -E_1,\ E_3\to E_5\to -E_3,\ E_4\to -E_6\to -E_4.$$

Then the Eisenstein form is given by:

\begin{equation}\label{}
\begin{cases}
h_{00} = \alpha E_{1} + {\sqrt{-1} \over 16}(\alpha + 63\beta)(E_{2} + E_{3} + E_{4} + E_{5}),\cr
h_{0} = \beta E_1 + {63 \sqrt{-1} \over 16}(\alpha - \beta) (E_{2} + E_{3} + E_{4} + E_{5}),\cr
h_{1} = {63 \sqrt{-1} \over 16}(\alpha - \beta) (E_{2} - E_{3} + E_{4} - E_{5}) + \beta E_6,\cr
h_{10} = {\sqrt{-1} \over 16} (\alpha - \beta)(E_{2} - E_{3} + E_{4} - E_{5}) + \alpha E_6,\cr
h_{1/2} = {56 \sqrt{-1} \over 16} (\alpha - 9\beta)(E_{2} + \sqrt{-1} E_{3} - E_{4} - \sqrt{-1} E_{5}),\cr
h_{3/2} = {72 \sqrt{-1} \over 16} (\alpha + 7\beta)(E_{2} - \sqrt{-1} E_{3} - E_{4} + \sqrt{-1} E_{5})\cr
\end{cases}
\end{equation}
where $\alpha, \beta$ are parameters.
We need the solution such that the constant terms of 
$h_{00}$, $h_0$, $h_{10}$, $h_{1}$, $h_{1/2}$, $h_{3/2}$ are $-1/2, 0,0,0,0,0$, respectively.  
Such Eisenstein form is given by:

\begin{equation}\label{eisen1}
\begin{cases}
h_{00} = \alpha E_{1} + {\sqrt{-1} \alpha \over 16} (E_{2} + E_{3} + E_{4}
+ E_{5}) = -{1 \over 2} + {2^{10} + 2 \over 61}q + \cdot \cdot \cdot ,\cr
h_{0} = {63 \sqrt{-1} \alpha \over 16} (E_{2} + E_{3} + E_{4} + E_{5})
= {4^{5}\cdot 63 \over 61}q + \cdot \cdot \cdot,\cr
h_{1} = {63 \sqrt{-1} \alpha \over 16} (E_{2} - E_{3} + E_{4} - E_{5})
= {2^{4}\cdot 63 \over 61} q^{1/2} + \cdot \cdot \cdot,\cr
h_{10} = {\sqrt{-1} \alpha \over 16} (E_{2} - E_{3} + E_{4} - E_{5})
+ \alpha E_{6} = {18 \over 61} q^{1/2} + \cdot \cdot \cdot,\cr
h_{1/2} = {56 \sqrt{-1} \alpha \over 16} (E_{2} + \sqrt{-1} E_{3} - E_{4} - 
\sqrt{-1} E_{5}) = {2^{4}\cdot 7^{2}\cdot 13 \over 61} q^{3/4} + \cdot \cdot \cdot,\cr
h_{3/2} = {72 \sqrt{-1} \alpha \over 16} (E_{2} - \sqrt{-1} E_{3} - E_{4} +
\sqrt{-1} E_{5})  = {18 \over 61} q^{1/4} + \cdot \cdot \cdot \cr
\end{cases}
\end{equation}
where 
$$\alpha = -{7!  \over (2\pi)^{7} B_{7}(1/4)} = 
{-4^{7}\cdot 6! \over (2\pi)^{7} 61}.$$

\subsection{Cusp forms}\label{cusp}

The obstruction space has dimension 4 and it
contains 2-dimensional subspace of cusp forms (Lemma \ref{dim}).
We shall calculate cusp forms in the obstruction space.
We consider the following two types:

\medskip

(A) $\{ \eta^{6}(\tau ) g_{\gamma} \}$;
\medskip

(B) $\{ \eta^{12}(\tau) g_{\gamma} \}$.

\subsection{Case (A)}\label{}
We denote by
$$F_{1} = G_{4}^{(0,1)}(\tau),\ F_{2} = G_{4}^{(1,0)}(\tau),\
F_{3} = G_{4}^{(1,1)}(\tau),$$
$$F_{4} = G_{4}^{(1,2)}(\tau),\ F_{5} = G_{4}^{(1,3)}(\tau),\
F_{6} = G_{4}^{(2,1)}(\tau)$$
the Eisenstein series of weight 4 and level 4.  Their Fourier expansion are as follows (Koblitz, \cite{Kob}, Chap. III, \S 3, Proposition 22):

\begin{equation}\label{}
\begin{cases}
F_{1} = c_{4}(1 - 2^{4} q^{2} + \cdot \cdot \cdot ); \cr
F_{2} = c_{4}(q^{1/4} + 2^{3} q^{2} + (1+3^{3}) q^{3/4} + 4^{3} q +
\cdot \cdot \cdot );\cr
F_{3} = c_{4}(\sqrt{-1} q^{1/4} - 2^{3} q^{2} - \sqrt{-1}(1+3^{3}) q^{3/4} + 4^{3} q +
\cdot \cdot \cdot );\cr
F_{4} = c_{4}(- q^{1/4} + 2^{3} q^{2} - (1+3^{3}) q^{3/4} + 4^{3} q +
\cdot \cdot \cdot );\cr
F_{5} = c_{4}(-\sqrt{-1} q^{1/4} - 2^{3} q^{2} + \sqrt{-1}(1+3^{3}) q^{3/4} + 4^{3} q +
\cdot \cdot \cdot );\cr
F_{6} = c_{4}(- 2^{4} q + \cdot \cdot \cdot ), \cr
\end{cases}
\end{equation}
where $c_{4} = {(2\pi)^{4} \over 2^{6}\cdot 4!}$.  
The acion of $S,T$ on $F_i$ is as follows: 

$$T(F_i)=F_i, i=1,6,\ T: F_2\to F_3\to F_4\to F_5\to F_2;$$
$$S: F_1 \to F_2\to F_1,\ F_3\to F_5\to F_3,\ F_4\to F_6\to F_4.$$

Recall that $\eta^6(\tau+1)=\sqrt{-1}\eta^6(\tau)$ and $\eta^6(-1/\tau)=\sqrt{-1}\tau^3\eta^6(\tau)$.
If we write $h_{\alpha} = \eta^{6}(\tau ) g_{\alpha}$, 
we need to find $\{ g_{\alpha}\}$ satisfying:

\begin{equation}\label{}
\begin{cases}
g_{00}(\tau + 1) = -\sqrt{-1}g_{00}(\tau), \ g_{0}(\tau + 1) = \sqrt{-1}g_{0}(\tau),\
g_{1}(\tau + 1) = \sqrt{-1}h_{1}(\tau),\cr
g_{10}(\tau + 1) = \sqrt{-1}g_{10}(\tau),\  
g_{1/2}(\tau + 1) = -g_{1/2}(\tau),\ g_{3/2}(\tau + 1) = g_{3/2}(\tau),\cr
g_{00}(-1/\tau) = {-\tau^{4} \over 16} (g_{00} + g_{0} + g_{1} + g_{10} + g_{1/2} + g_{3/2}),\cr
g_{0}(-1/\tau) = {- \tau^{4} \over 16} (63 g_{00} - g_{0} - g_{1} + 63 g_{10} - 9 g_{1/2} + 7 g_{3/2}),\cr
g_{1}(-1/\tau) = {- \tau^{4} \over 16} (63 g_{00} - g_{0} - g_{1} + 63 g_{10} + 9 g_{1/2} - 7 g_{3/2}),\cr
g_{10}(-1/\tau) = {- \tau^{4} \over 16} (g_{00} + g_{0} + g_{1} + g_{10} - g_{1/2} - g_{3/2}),\cr
g_{1/2}(-1/\tau) = {- \tau^{4} \over 16} (56 g_{00} - 8 g_{0} + 8g_{1} - 56 g_{10}),\cr
g_{3/2}(-1/\tau) = {- \tau^{4} \over 16} (72 g_{00} + 8 g_{0} - 8 g_{1} - 72 g_{10}).\cr
\end{cases}
\end{equation}

Then one dimensional subspace of cusp forms is given by

\begin{equation}\label{cuspA}
\begin{cases}
h_{00} = - {\beta \eta(\tau)^{6} \over 16} (F_{2} + \sqrt{-1} F_{3} - F_{4} - \sqrt{-1} F_{5}) = 
-{\beta c_{4} \over 16}(4(3^{3} + 1)q + \cdot \cdot \cdot),\cr
h_{0} = -{7 \beta \eta(\tau)^{6} \over 16} (F_{2} + \sqrt{-1} F_{3} - F_{4} - \sqrt{-1} F_{5}) = 
-{7\beta c_{4} \over 16}(4(3^{3}+1)q + \cdot \cdot \cdot),\cr
h_{1} = {7 \beta \eta(\tau)^{6} \over 16} (F_{2} - \sqrt{-1} F_{3} - F_{4} + \sqrt{-1} F_{5}) = 
{7\beta c_{4} \over 16}(4q^{1/2} + \cdot \cdot \cdot ),\cr
h_{10} = {\beta \eta(\tau)^{6} \over 16} (F_{2} - \sqrt{-1} F_{3} - F_{4} + \sqrt{-1} F_{5}) = 
{\beta c_{4} \over 16}(4q^{1/2} + \cdot \cdot \cdot ),\cr
h_{1/2} = 0,\cr
h_{3/2} = \beta \eta(\tau)^{6} (F_{1} - F_{6})
= \beta c_{4}(q^{1/4} + \cdot \cdot \cdot )\cr
\end{cases}
\end{equation}
where $\beta$ is a parameter.

\subsection{Case (B)}\label{}
Let 
$$G_{1} = G_{1}^{(0,1)}(\tau), \ G_{2} = G_{1}^{(1,0)}(\tau),\
G_{3} = G_{1}^{(1,1)}(\tau),$$
$$G_{4} = G_{1}^{(1,2)}(\tau),\ G_{5} = G_{1}^{(1,3)}(\tau),\
G_{6} = G_{1}^{(2,1)}(\tau)$$
be the Eisenstein series of weight 1 and level 4. 

The action of $S, T$ on $G_{i}$ is:
$$T(G_{1}) = G_{1}, G_{2} \rightarrow G_{3} \rightarrow G_{4}
\rightarrow G_{5} \rightarrow G_{2}, T(G_{6}) = -G_{6};$$
$$S : G_{1} \rightarrow G_{2} \rightarrow -G_{1}, 
G_{3} \rightarrow G_{5} \rightarrow -G_{3},
G_{6} \rightarrow G_{4} \rightarrow -G_{6}.$$
Write $h_{\alpha} = \eta^{12}(\tau)g_{\alpha}$.

Since $\eta^{12}(\tau + 1) = -\eta^{12}(\tau), \eta^{12}(-1/\tau)
= - \tau^{6} \eta^{12}(\tau)$, $\{ g_{\alpha}\}$ satisfies the following:

\begin{equation}\label{}
\begin{cases}
g_{00}(\tau + 1) = -g_{00}(\tau),\ g_{0}(\tau + 1) = -g_{0}(\tau),\
g_{1}(\tau + 1) = g_{1}(\tau),\cr
g_{10}(\tau + 1) = g_{10}(\tau),\ g_{1/2}(\tau + 1) = \sqrt{-1} g_{1/2}(\tau),\
g_{3/2}(\tau + 1) = -\sqrt{-1} g_{3/2}(\tau),\cr
g_{00}(-1/\tau) = {\sqrt{-1} \tau \over 16} (g_{00} + g_{0} + g_{1} +
g_{10} + g_{1/2} + g_{3/2}),\cr
g_{0}(-1/\tau) = {\sqrt{-1} \tau \over 16} (63 g_{00} - g_{0} - g_{1} +
63 g_{10} - 9 g_{1/2} + 7 g_{3/2}),\cr
g_{1}(-1/\tau) = {\sqrt{-1} \tau \over 16} (63 g_{00} - g_{0} - g_{1} +
63 g_{10} + 9 g_{1/2} - 7 g_{3/2}),\cr
g_{10}(-1/\tau) = {\sqrt{-1} \tau \over 16} (g_{00} + g_{0} + g_{1} +
g_{10} - g_{1/2} - g_{3/2}),\cr
g_{1/2}(-1/\tau) = {\sqrt{-1} \tau \over 16} (56 g_{00} - 8 g_{0} + 
8 g_{1} - 56 g_{10}),\cr
g_{3/2}(-1/\tau) = {\sqrt{-1} \tau \over 16} (72 g_{00} + 8 g_{0} - 
8 g_{1} - 72 g_{10}).\cr
\end{cases}
\end{equation}

$\{ G_{i} \}$ are not linearly independent, but there are no cusp forms
of weight 1 and level 4.  Hence $g_{\alpha}$ can be written as a linear
combination of $\{ G_{i} \}$.
Then solutions are

\begin{equation}\label{}
\begin{cases}
g_{00} = - {\sqrt{-1}(\alpha + \beta) \over 16} 
(G_{2} - G_{3} + G_{4} - G_{5}) + \alpha G_{6},\cr
g_{0} = - {\sqrt{-1}(63\alpha - \beta) \over 16} 
(G_{2} - G_{3} + G_{4} - G_{5}) + \beta G_{6},\cr
g_{1} = \beta G_{1} - {\sqrt{-1}(63\alpha - \beta) \over 16} 
(G_{2} + G_{3} + G_{4} + G_{5}),\cr
g_{10} = \alpha G_{1} - {\sqrt{-1}(\alpha +\beta) \over 16} 
(G_{2} + G_{3} + G_{4} + G_{5}),\cr
g_{1/2} = {\sqrt{-1}(56 \alpha - 8 \beta) \over 16}
 (G_{2} - \sqrt{-1} G_{3} - G_{4} + \sqrt{-1} G_{5}),\cr
g_{3/2} = {\sqrt{-1}(72 \alpha + 8 \beta) \over 16} 
(G_{2} + \sqrt{-1} G_{3} - G_{4} - \sqrt{-1} G_{5}).\cr
\end{cases}
\end{equation}
The Fourier coefficients of $G_{i}$ are as follows
(Schoenberg \cite{S}, Chap. VII, \S 2):

\begin{equation}\label{}
\begin{cases}
G_{1} = a_{0} + \pi q + \cdot \cdot \cdot \ , \cr
G_{2} = b_{0} - {\pi \sqrt{-1} \over 2}(q^{1/4} + q^{1/2} + q + \cdot \cdot), \cr
G_{3} = b_{0} - {\pi \sqrt{-1} \over 2}(\sqrt{-1}q^{1/4} - q^{1/2} + q + \cdot \cdot), \cr
G_{4} = b_{0} - {\pi \sqrt{-1} \over 2}(-q^{1/4} + q^{1/2} + q + \cdot \cdot), \cr
G_{5} = b_{0} - {\pi \sqrt{-1} \over 2}(-\sqrt{-1}q^{1/4} - q^{1/2} + q + \cdot \cdot), \cr
G_{6} = c_{0} + \pi q^{1/2} + \cdot \cdot \cdot \ , \cr
\end{cases}
\end{equation}
where
$$a_{0} = {1 \over 4} \lim_{s \to 0} [\zeta(1+s,1/4) - \zeta(1+s,3/4)],$$
$$b_{0} = {-\pi \sqrt{-1} \over 4}[\zeta(0,1/4) - \zeta(0,3/4)],$$
$$c_{0} = {-\pi \sqrt{-1} \over 4}[\zeta(0,1/2) - \zeta(0,1/2)] = 0.$$
Since $\zeta(0,a) = 1/2 - a$ (Whittaker, Watson \cite{W},
13.21), $b_{0} = -\pi \sqrt{-1}/8$.
On the other hand, 
$$\lim_{s \to 0} [\zeta(s+1, a) - 1/s] = -{\Gamma'(a) \over \Gamma(a)}$$
(Whittaker, Watson, \cite{W}, 13.21).  Hence
$$a_{0} = {\Gamma(1/4) \Gamma'(3/4) - 
\Gamma'(1/4) \Gamma(3/4) \over 4\Gamma(1/4) \Gamma(3/4)}.$$
By $\Gamma(z) \Gamma(1-z) = \pi / \sin \pi z$ (Whittaker, Watson \cite{W}, 12.14), we have
$a_{0} = \pi/4$.
Thus we have

\begin{equation}\label{cuspB}
\begin{cases}
h_{00}=\eta(\tau)^{12} g_{00} = {\pi (7\alpha - \beta) \over 8} q +
\cdot \cdot \cdot\ ,\cr
h_0=\eta(\tau)^{12} g_{0} = - {9 \pi (7\alpha - \beta) \over 8} 
q + \cdot \cdot \cdot \ ,\cr
h_1=\eta(\tau)^{12} g_{1} = -{9 \pi (7\alpha - \beta) \over 32}
q^{1/2} + \cdot \cdot \cdot  \ ,\cr
h_{10}=\eta(\tau)^{12} g_{10} = {\pi (7 \alpha - \beta) \over 32} 
q^{1/2} + \cdot \cdot \cdot \ ,\cr
h_{1/2}=\eta (\tau)^{12} g_{1/2} = \pi (7 \alpha - \beta) q^{3/4} + \cdot
\cdot \cdot \ , \cr
h_{3/2}=\eta(\tau)^{12} g_{3/2} = 0\cdot q^{1/4} + \cdot \cdot\cdot \ , \cr
\end{cases}
\end{equation}

\subsection{Theorem}\label{multiplicative}
{\it Let}
$$D = m_1\calD_1 + m_{10}\calD_{10} + m_{1/2}\calD_{1/2} + m_{3/2}\calD_{3/2}.$$
{\it Then $D$ is a divisor of meromorphic automorphic form $F$ on $\calD$ if
$$7m_1 + m_{10} + 4m_{3/2} = 0, \quad -9m_1 + m_{10} + 32m_{1/2} = 0,$$
that is
$$m_{10} = 9m_1 -32m_{1/2}, \quad m_{3/2} = 8m_{1/2} - 4m_1.$$
In this case, the weight of $F$ is given by}
$$2\cdot 3^2 m_1 + 2^5\cdot 5 m_{1/2}.$$

\begin{proof}
The first assertion follows from Theorem \ref{freitag} and the equations \ref{cuspA}, \ref{cuspB}.  By using Theorem \ref{freitag} and the equation \ref{eisen1}, we can see that 
the weight of $F$ is given by
$${2^4\cdot 63 \over 61}m_1 + {18 \over 61}m_{10} + {2^4\cdot 7^2 \cdot 13 \over 61}m_{1/2} + {18\over 61}m_{3/2} = 2\cdot 3^2 m_1 + 2^5\cdot 5 m_{1/2}.$$
\end{proof}

\subsection{Corollary}\label{multi}
{\it There exists a meromorphic
 automorphic forms $\Phi$ on $\calD$ of weight $18$ whose divisor is}  
$${\calD}_{1} + 9{\calD}_{10} - 4{\calD}_{3/2}.$$

\begin{proof}
Put $m_1 =1, m_{1/2}=0$ in Theorem \ref{multiplicative}.
\end{proof}
\subsection{Theorem}\label{main2}
{\it Let $I$ be a maximal totally isotropic subspace and let $V =\la I, \kappa\ra$ as in Definition {\rm \ref{invariant}}.  Let
$F_V$ be the aditive lifting corresponding to $V$.  Then the divisor of $F_V$ is}
$$ \sum_{a\in V,\ q(a)=1} \calD_a - \sum_{c\in I}(2\calD_{c+\alpha} + 2\calD_{c+\alpha+\kappa}).$$ 

\begin{proof}
Consider the product $\Psi$ of all $135$ $F_V$ which has weight $135\cdot 2$.
It follows from Theorem \ref{main1} that $\Psi$ vanishes
along $\calD_1$ with at least multiplicity $15$ and along $\calD_{10}$ with at least 
multiplicity $135$ and has poles along $\calD_{3/2}$ with multiplicity $60$.
Then the quotient $\Psi / \Phi^{15}$ is a holomorphic automorphic form of weight zero (Corollary \ref{multi}).  
Now the Koeher principle implies that $\Psi/\Phi^{15}$ is constant.
\end{proof}
\subsection{Corollary}\label{cubicrelation}
{\it The $135$ $F_V$ satisfies $63$ cubic relations corresponding to $63$ non-zero isotropic vectors.}

\begin{proof}
We use the notation as in  Lemma \ref{cubic}.  
By applying Theorem \ref{main2}, the automorphic forms $F_{V_1}F_{V_2}F_{V_3}$ and
$F_{V_4}F_{V_5}F_{V_6}$ have the same divisor and hence
$$F_{V_1}F_{V_2}F_{V_3}= c F_{V_4}F_{V_5}F_{V_6}\  (c \in \bbC).$$
\end{proof}

\subsection{Corollary}\label{restriction}
{\it The restriction of $F_V$ to $\calB$ has the divisor 
$$2\calB_{\kappa} + \sum_{a\in V,\ q(a)=1, a\not=\kappa} 2\calB_a - \sum_{c\in I} 4\calB_{c+\alpha}.$$
}

\begin{proof}
Note that $\calD_{c+\alpha} \cap \calB = \calD_{c+\alpha+\kappa} \cap \calB = \calB_{c+\alpha}$ and the fact  (\ref{fact1}).
Then the assertion follows.
\end{proof}

By definition, $F_V$ is invariant with respect to $\tilde{\Gamma}$ but not invariant under $\tilde{\Gamma}\cdot \bbZ/2\bbZ$.
However $F_V(\rho(\omega)) = F_V(\sqrt{-1}\omega) = (\sqrt{-1})^2F_V(\omega) = -F_V(\omega)$ for any $V$.  Hence
for totally isotropic subspaces $I, I'$, the quotient $F_V/F_{V'}$ $(V'=\la I', \kappa\ra)$ is invariant with respect to 
$\tilde{\Gamma}\cdot \bbZ/2\bbZ$.
Thus we can determine the divisor of $F_V/F_{V'}$ on $\calB/(\tilde{\Gamma}\cdot \bbZ/2\bbZ)$. 
Let
$$2\calB_{\kappa} + \sum_{a'\in V',\ q(a')=1, a'\not=\kappa} 2\calB_{a'} - \sum_{c'\in I'} 4\calB_{c'+\alpha'}$$
be the divisor of $F_{V'}$ on $\calB$.
Let $\pi : \calB \to \calB/(\tilde{\Gamma}\cdot \bbZ/2\bbZ)$ be the natural projection.  Then $\pi$ is branched over $\calH_{\alpha}$
and its branch degree is two (resp. four) if $\alpha$ is of type $(1)$ (resp. type $(3/2)$) (Lemma \ref{branch}).  
Hence the vanishing order or pole order of $F_V/F_V'$ along
$\calH_{\alpha}$ is given by dividing  the vanishing order or pole order along $\calB_{\alpha}$ by its branch degree. 
Thus we have:

\subsection{Corollary}\label{main3}
{\it The divisor of $F_V/F_{V'}$ on $\calB /(\tilde{\Gamma}\cdot \bbZ/2\bbZ)$ is given by
$$\sum_{a\in V,\ q(a)=1, a\not=\kappa} \calH_a - \sum_{a'\in V',\ q(a')=1, a'\not=\kappa} \calH_{a'} + 
\sum_{c'\in I'} \calH_{c'+\alpha'} - \sum_{c\in I} \calH_{c+\alpha}.$$
}

\medskip

\section{Automorphic forms and G{\" o}pel invariants}

The purpose of this section is to prove the following:

\subsection{Theorem}\label{main4}
{\it The linear system defined by $135$ meromorphic automorphic forms $F_V$ gives a 
rational map
$$\varphi : \calB/(\tilde{\Gamma}\cdot \bbZ/2\bbZ) --\to \bbP^{14}$$
which is birational onto its image.  The image of $\varphi$ satisfies $63$ cubic relations.
Under the identification between $P_2^7$ and $\calB/(\tilde{\Gamma}\cdot \bbZ/2\bbZ)$
by a birational isomorphism, the map $\varphi$ coincides with $\psi$ given by G{\" o}pel invariants in Theorem} \ref{Goepel}.

\begin{proof}
We shall show that the map $\varphi$ coincides with $\psi$.  Then the remaining assertions 
follow from Theorem \ref{Coble}, Corollary \ref{cubicrelation}. 

Under the isomorphism $\hat{p}$ in Proposition \ref{kondo1}, we identify $(\hat{P}_2^7)_1$ with 
$(\calB\setminus \calH_h)/(\tilde{\Gamma}\cdot \bbZ/2\bbZ)$ up to subvarieties of codimension 2, 
and then we shall compare $\psi$ and $\varphi$.

Recall that the quadratic space $E_7/2E_7$ is canonically isomorphic to the subspace $A_{L_-}'$ of $q_{L_-}$ 
(Lemma \ref{e7}).  
Hence G\"opel subspaces bijectively correspond to totally isotropic subspaces in $q_{L_-}$.
Consider a totally isotropic subspace $A$ of dimension 2 in $q_{L_-}$.  Then there exist exactly three totally isotropic subspaces
$I_1, I_2, I_3$ in $q_{L_-}$ containing $A$.  Let $M_1, M_2, M_3$ be the corresponding G\"opel subspaces.
For example, $M_1, M_2, M_3$ are given by (\ref{goepel-3}), (\ref{goepel-4}), (\ref{goepel-2}) respectively.
Put $V_i= \la I_i, \kappa \ra$.  Then it follows from Lemma \ref{blowup} and Corollary \ref{main3} that
$(\hat{G}_{M_1}/\hat{G}_{M_2}) = (F_{V_1}/F_{V_2})$ as divisors.  
Let $\alpha$ be a positive root in $M_3$ not appeared in $M_1$ and $M_2$.
We can see that $\hat{G}_{M_1}/\hat{G}_{M_2}$ takes value 1 on the divisor $D_{\alpha}$.  
We denote by $\calB_a$ a component of the divisor
of  $F_{V_3}$ corresponding to $D_{\alpha}$.
We may assume that
$G_{M_1} - G_{M_2} = G_{M_3}$ (see Proposition \ref{Goepel}, (3)).
Then  $\theta_{V_1} - \theta_{V_2} = \pm \theta_{V_3}$ (see Lemma \ref{linear-rel}) and hence 
$F_{V_1} - F_{V_2}$ vanishes on the divisor of $F_{V_3}$.
This implies that $F_{V_1}/F_{V_2}$ takes value 1 on the divisor $\calH_a$.  Hence $\hat{G}_{M_1}/\hat{G}_{M_2} = F_{V_1}/F_{V_2}$.
For any pair of G\"opel subspaces $M_1', \ M_2'$, the ratio $\hat{G}_{M_1'}/\hat{G}_{M_2'}$ can be written as the product of
$\hat{G}_{M_1}/\hat{G}_{M_2}$ as above type.  Thus we have finished the proof of Theorem \ref{main4}.
\end{proof}

\subsection{Remark}\label{}
In the above theorem, we do not consider the hyperelliptic locus.  The 36 subvarieties of $P_2^7$ of codimension 2 
in Lemma \ref{28/36} bijectively correspond to 36 divisors $\calH_{\alpha}$ ($\alpha$ is type $(3/2)$) (see Remark \ref{e7-5}).  
Let $g : \tilde{P}_2^7 \to \hat{P}_2^7$ be the blow up along 36 subvarieties.  Then we can see that 
$g^*(\tilde{G}_{V_1}/\tilde{G}_{V_2})$ 
has the same divisor as $F_{V_1}/F_{V_2}$ by using Lemma \ref{28/36}.  The author conjectures that the map $\hat{p}$ can be extended to the hyperelliptic locus.


\end{document}